\pdfoutput=1 % added by arXiv admin
\documentclass[a4paper]{amsart}

\usepackage{amsfonts,amsmath,amsthm}
\usepackage{mathtools}

\usepackage{graphicx}

\title[The iterated refinement Newton method]
{Newton's method in practice II: \\ 
The iterated refinement Newton method and 
near-optimal complexity for 
finding all roots of some polynomials of very large degrees}
\author{Marvin Randig}
\author{Dierk Schleicher}
\author{Robin Stoll}

\theoremstyle{definition}

\newcommand{\eps}{\varepsilon}

\newcommand{\id}{\operatorname{id}}

\newcommand{\C}{{\mathbb C}}

\newcommand{\disk}{{\mathbb D}}

\renewcommand{\phi}{\varphi}

\newcommand{\reminder}[1]{\textsl{#1}}%\protect\marginpar{$\star$}}

\newcommand{\hide}[1]{}

\begin{document}

\begin{abstract}
We present a practical implementation based on Newton's method to find all roots of several families of complex polynomials of degrees exceeding one billion ($10^9$) so that the observed complexity to find all roots is between $O(d\ln d)$ and $O(d\ln^3 d)$ (measuring complexity in terms of number of Newton iterations or computing time). All computations were performed successfully on standard desktop computers built between 2007 and 2012.
\end{abstract}

\maketitle

\tableofcontents

\section{Introduction}

It has been known since Gauss that every complex polynomial $p$ in a single complex variable splits into linear factors, and since Ruffini and Abel that there is no method based on finite concatenation of $n$-th roots to find these factors algebraically in general. Therefore, numerical approximation methods are required to find the roots of polynomials. 

\begin{figure}[htbp]
\includegraphics[width=0.6\textwidth]{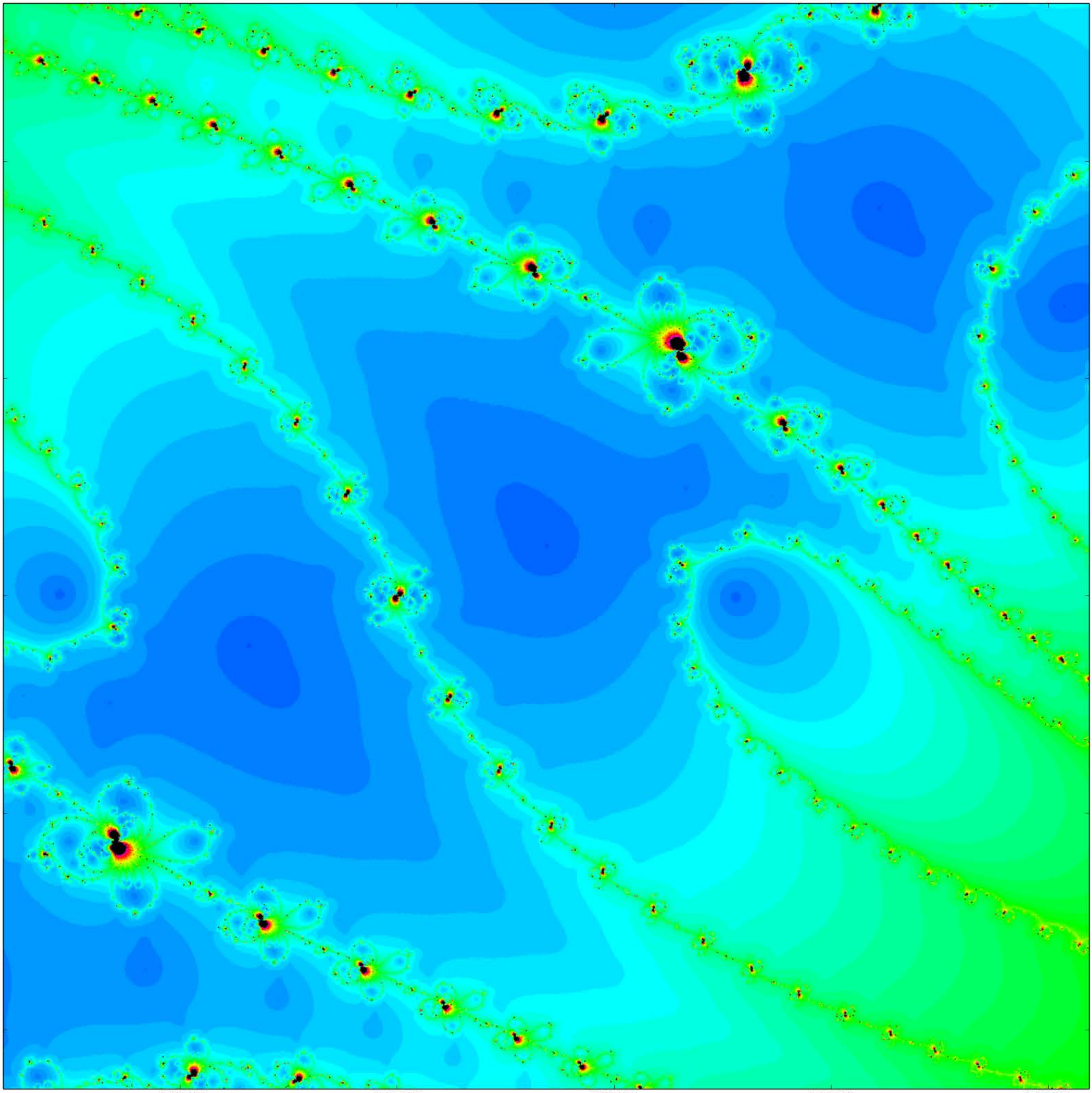}
\caption{The Newton dynamics of a polynomial of degree $2^{27}>10^8$ (detail). The colors illustrate the number of iterations until a root is found to high precision; adjacent color contour lines show a single iteration. The scale is between 50 iterations (red) and 0 iterations (dark blue); bright green stands for 25 iterations. Compare with Figure~\ref{Fig:NewtonMarvin_1_24}. }
\label{Fig:NewtonMarvin_1_19}
\end{figure}

Surprisingly, even today it is not clear how to find all roots most efficiently. There is a theoretically best possible algorithm due to Pan \cite{Pan02}, but it is not applicable in practice. Perhaps the most efficient practical implementation is \texttt{MPSolve 3.0} by Bini and Robol \cite{BiniRobol} based on the Ehrlich-Aberth iteration, but it is still lacking good theory and especially a proof that it is generally convergent.

Newton's method is very good at finding roots locally, i.e.\ once good approximations are known, but it has a reputation of being difficult to predict globally, for instance because of the ``chaotic'' nature of the iteration, and because of the possibility of having open sets of starting points that converge to attracting periodic orbits, rather than to roots.

The situation has improved substantially in recent times. Indeed, a pleasant feature of Newton's method is that one does not have to choose between theory and practical implementations: meanwhile it has good theory (see \cite{HSS,BLS} on where to start the iteration, and \cite{NewtonEfficient} with refinements in \cite{BAS} on estimates on its complexity that are near-optimal under certain assumptions), and it can be implemented quite successfully in practice for polynomials of very high degrees. The latter aspect is the key theme of the present manuscript.

In \cite{NewtonRobin1}, we showed that Newton's method can be used to find all roots of various complex polynomials efficiently, even up to degrees greater than one million (where the polynomials were selected only on the criterion that they could be evaluated efficiently). In this note, we demonstrate that the computational complexity (measured in computing time) can be near-optimal for a wide range of degrees: all roots can be found in computing time between $O(d\ln d)$ and $O(d\ln^3d)$, and in a number of Newton iterations between $O(d\ln d)$ and $O(d\ln^2 d)$, with small constants, so that all roots can be found for degrees up to $2^{30}$ (greater than one billion), and even higher degrees seem feasible. 

We present an algorithm based on Newton's method that turns out to find all roots, the families of polynomials studied, and we describe the outcome of computer experiments. We also explain how to verify that indeed all roots of our polynomials have been found.

This note should be seen as a continuation and improvement on \cite{NewtonRobin1}, not as a survey on root finding for polynomials. There are interesting references on the latter, in particular the survey articles by McNamee~\cite{McNamee,McNameeBook}, as well as studies by Pan \cite{Pan,Pan02}, Renegar~\cite{Renegar}, and especially MPSolve 3.0 from \cite{BiniRobol} and Eigensolve from \cite{Fortune}.

\textbf{Acknowledgements}. We are very pleased to be able to report support, encouragement, and helpful suggestions from numerous friends and colleagues, especially Dario Bini, Marcel Oliver, Victor Pan, Simon Schmitt, and Michael Stoll, and we are most grateful to them all. 

\section{Background on Newton dynamics}

There are a number of theoretical results available on the global dynamics of Newton's method $N_p=\id-p/p'$ for a given polynomial $p$. All require that a disk is known that contains all roots of $p$. After appropriate rescaling, we may assume that all roots of $p$ are contained in $\disk$. 

We showed in \cite{HSS} that there is a universal set of $1.1d\ln^2d$ starting points, placed on $O(\ln d)$ concentric circles containing $O(d\ln d)$ points each, so that Newton's method started at these points will find all roots of all polynomials of degree~$d$. This was improved in \cite{BLS} to a probabilistic set of starting points containing only $O(d(\ln \ln d)^2)$ starting points that find all roots with arbitrarily high probability. Our earlier experiments in \cite{NewtonRobin1} started with $4d$ or $8d$ equidistributed points on a single circle; in many cases, $4d$ points were sufficient, but always $8d$ points were enough: this shows that in practice $O(d)$ starting points seem to suffice, even though this is not supported by theory (and there may well be special polynomials that require more starting points). Of course, any additional factors like $\ln^2 d$ are hard to detect numerically, even when the degree ranges up to $2^{30}$.

Theoretical upper bounds on the expected number of required Newton iterations were given in \cite{NewtonEfficient,BAS}; these are, up to polylogarithmic factors, on the order of $d^3$ or even $d^2$.

However, these previous methods all require at least $\Omega(d^2)$ iterations in order to find all roots: the methods rely on controlled Newton dynamics away from $\disk$, so all starting points $z$ are uniformly bounded away from the disk: $|z|>r>1$ for a uniform $r$. However, on this domain we essentially have $N_p(z)\approx dz/(d-1)$ \cite[Lemma~3]{HSS}, so each orbit needs at least $\Omega(d\ln r)$ iterations until it even enters $\disk$, and since at least $d$ orbits are required to find $d$ roots, we obtain a lower bound for the complexity of $\Omega(d^2)$ Newton iterations. The upper bounds in \cite{BAS} are thus best possible, again up to polylogarithmic factors. 

One may try to remedy this problem by starting the iteration on a circle of radius, say, $1+1/d$. Such an approach comes with various problems. One is that we are losing the theory to guarantee that all roots will be found. Another one is that this will be helpful only if a very tight bound for the smallest disk containing all roots is known, and only if most roots are near the boundary of that disk. This is not the case for instance for all polynomials considered by us (but it is true if the coefficients of the polynomial are independently randomly distributed \cite{ErdoesTuran,Arnold}).

\section{The Iterated Refinement Newton Method}
\label{Sec:IteratedRefinementNewton}

\begin{figure}[p]
\includegraphics[width=0.65\textwidth]{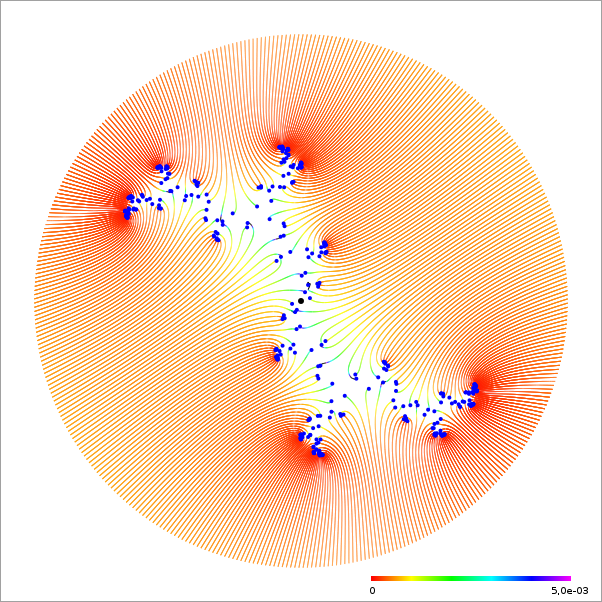}
\caption{The Newton iteration for $400$ starting points on a circle of radius $r=2$ (here for a polynomial of degree $4096$, so we do not have enough starting points to find all roots; the polynomial shown here describes periodic points of period dividing $12$ of $z\mapsto z^2+i$). The apparent lines connect orbits under the Newton dynamics; colors indicate the number of iterations until an approximate root is found.
The behavior of the iterations outside of the disk containing all roots is very parallel and ``wasteful'', but required to carry over the control from the circle of starting points to the interesting dynamics on $\disk$. Observe that even if we had a very precise bound on the smallest disk containing all roots, this would not help much as most roots are away from the boundary of this disk.  }
\label{Fig:ParallelIterations}

\includegraphics[width=0.65\textwidth]{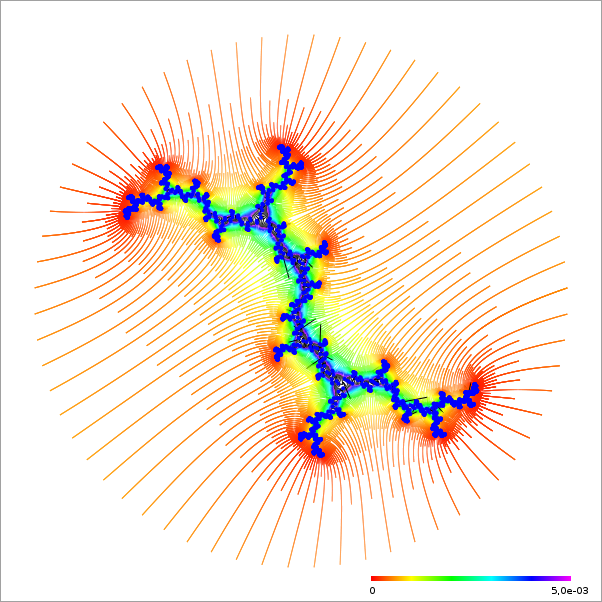}
\caption{The Iterated Refinement Newton Algorithm: we start with $64$ points and refine only when adjacent orbits stop moving in parallel.  As a result, most Newton iterations are spent near the roots, and far fewer iterations are required. The polynomial is the same as in Figure~\ref{Fig:ParallelIterations} (degree $4096$). }
\label{Fig:NewtonRefinement}
\end{figure}

One key observation is that the lower bound of the number of iterations comes from the iterations outside of $\disk$, that is before the interesting dynamics even starts, and where all starting points are on very similar ``parallel'' orbits --- a necessary consequence of the fact that we start the Newton dynamics at points with controlled dynamics. Typical initial orbits of the Newton iteration are shown in Figure~\ref{Fig:ParallelIterations}. The approach we take in the experiments described in this paper is that on domains where orbits are ``parallel'', fewer orbits are required, as indicated in Figure~\ref{Fig:NewtonRefinement}.

More precisely, we start the Newton dynamics with a fixed number of, say, $64$ Newton orbits on a circle away from $\disk$ where we have good control. Along each triple of adjacent orbits, we compare the shapes of the triangles formed by the three points on adjacent orbits. As long as this shape remains nearly constant, so that the three orbits form nearly similar triangles, we keep iterating these three orbits, assuming that they represent similar dynamics for all orbits that might run between these three ``representative'' orbits. Once the triangles start to deform, we split up these orbits by inserting an additional orbit half way between each of the two pairs of adjacent orbits; see Figure~\ref{Fig:NewtonRefinement}. 

Of course, this requires a heuristic refinement threshold parameter that describes the threshold of deformation after which parallel orbits are refined. We use the cross-ratio between three adjacent orbits with respect to $\infty$: that is, 
\begin{equation}
t_{n,i}= \frac{z_{n,i-1}-z_{n,i}}{z_{n,i+1}-z_{n,i}}
\end{equation}
where the index $n$ counts iterations and the index $i$ counts the circular order of the orbits (with due modifications at the end to turn the linear order of natural numbers into a circular order). That is, without refinement, we have $z_{n+1,i}=N_p(z_{n,i})$; and when refinements occur, these are inserted into the circular order (and the indices $i$ are shifted accordingly). 

Our refinement condition is essentially based on $\alpha_{n,i}:=|\ln(t_{n,i}/t_{n_0,i})|$, where $n_0\le n$ is the most recent iteration when a refinement occurred for this orbit (properly accounting for index shifts because of refinements that might have occurred for other orbits).  

The algorithm is determined by a number of parameters that we determined heuristically:
\begin{itemize}
\item
the number $N_0$ of initial orbits (we usually use $N_0=64$, but this should not make much of a difference);
\item
the maximal number of iterations a particular orbit can run for (we use $10d$); more important than this number is a good detection of periodic cycles (if a positive fraction of orbits runs to the maximal number of allowed iterations, we end up at complexity $O(d^2)$); 
\item
the maximal number of refinement generations (we use $\lg_2(4d/N_0)$; that means that if all refinements have taken place, we end up with $4d$ orbits);
\item
the refinement threshold $R$: that is, the maximal value of $\alpha_{n,i}$ that is allowed before further refinement occurs. Most of the time, we use the value of $R=0.05$, but for some experiments we had to increase sensitivity to $R=0.0005$;  
\item
and a stopping criterion when success is declared that some orbit found a root: we stop when $|N_p(z)-z|=|p(z)/p'(z)|<\eps_{\text{stop}}$; in practice, we use $\eps_\text{stop}=10^{-15}$ or sometimes $10^{-16}$ (while we used the 
\texttt{long double} data type with numerical precision of about $18$ relevant digits). In experiments for very high periods, we used higher precision arithmetic and chose $\varepsilon_\text{stop} = 10^{-18}$. 
\end{itemize}

%\newpage

\looseness-1
Finally, one needs a post-processing step: many roots will be found by several different orbits, and we have to make sure these roots are accounted for only once. Here we use very simple-minded heuristics and declare that two orbits found different roots if they terminate at distance greater than $\eps_\text{root}$; here we use $\eps_\text{root}=10^{-14}$ (sometimes $\eps_\text{root}=10^{-15}$ or $\eps_\text{root}=10^{-16}$). Of course, this requires an efficient way to sort the roots by mutual proximity (when all roots are found by an algorithm that is linear up to log-factors, then for naive sorting procedures it is easy to spend more time on sorting than on root finding). 

Of course all quantities are heuristic, and especially the choice of stopping criterion and post-processing (distinguishing different roots found) have obvious difficulties in the presence of high degrees or near-multiple roots. Our point is not that this is a strong part of the algorithm, but that even with such a simple-minded approach all roots can be found for all the polynomials investigated, and for extremely large degrees. An improved approach for more ``challenging'' families of polynomials, especially with near-multiple roots, is the subject of ongoing research that will be discussed elsewhere.

In Section~\ref{Sec:AllRootsFound} we describe several methods that establish a posteriori that all roots have been found with high accuracy.

\goodbreak

%\newpage

\section{The polynomials investigated and the results of  experiments}

We investigate (a subset of) the same polynomials we already discussed in \cite{NewtonRobin1}. We will briefly introduce the relevant families of polynomials and then present the results of the experiments. These polynomials have been chosen solely for the purpose of fast evaluation (in terms of recursion): our focus is on root finding not efficient polynomial evaluation, so we chose polynomials where evaluation is easy. Therefore, comparing our experiments with other experiments on different polynomials in terms of computing time yields an unfair bias in our favor, but comparisons in terms of the required number of Newton iterations should be meaningful.

\subsection{Centers of hyperbolic components of the Mandelbrot set}

These polynomials in a variable $c$ are defined by recursion $p_0=0$, $p_{n+1}=p_n^2+c$, so we have $p_1(c)=c$, $p_2(c)=c^2+c$, $p_3(c)=(c^2+c)^2+c$ etc., so that $p_n$ has degree $2^{n-1}$. Roots of $p_n$ are those parameters $c$ for which the iteration $z_0=0$, $z_{n+1}=z_n^2+c$ is periodic with period (dividing) $n$; this implies that $p_k$ divides $p_n$ when $k$ divides $n$ (using the fact that all roots of all $p_n$ are simple; compare \cite[Section~19]{Orsay}). 

\begin{figure}[htbp]
\framebox{
\includegraphics[width=\textwidth,trim=8 0 5 0,clip]{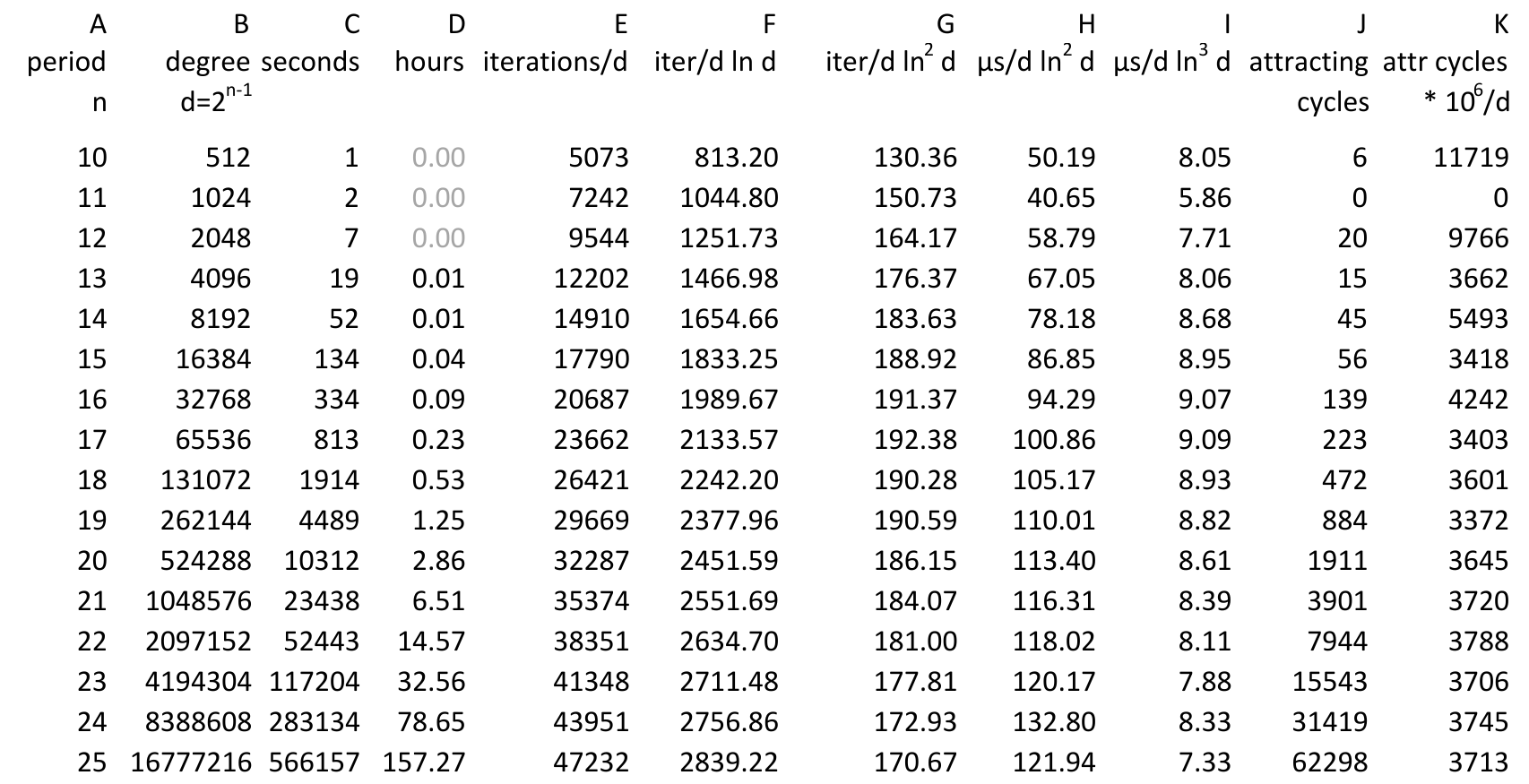} 
}

\caption{Results for finding all centers of hyperbolic components of the Mandelbrot set for period $n\le 25$, i.e.\ for degrees up to $2^{24}>1.6\cdot 10^7$. The first four columns specify period $n$ and degree $d=2^{n-1}$, as well as total computing time in seconds and hours. Column E gives the required number of Newton iterations (divided by $d$), while columns F and G scale this with respect to $\ln d$ and $\ln^2d$. Columns H and I specify computing time (in microseconds) divided by $d\ln ^2d$ and $d\ln ^3d$. Column J shows the number of orbits that converge to attracting cycles, and the final column K shows the same number times $10^6/d$. Refinement threshold $R=0.0005$, all roots were found. }
\label{Fig:MandelbrotComplexityTable}
\end{figure}

\begin{figure}[htbp]
\includegraphics[trim=55 530 88 75,clip,width=.8\textwidth]{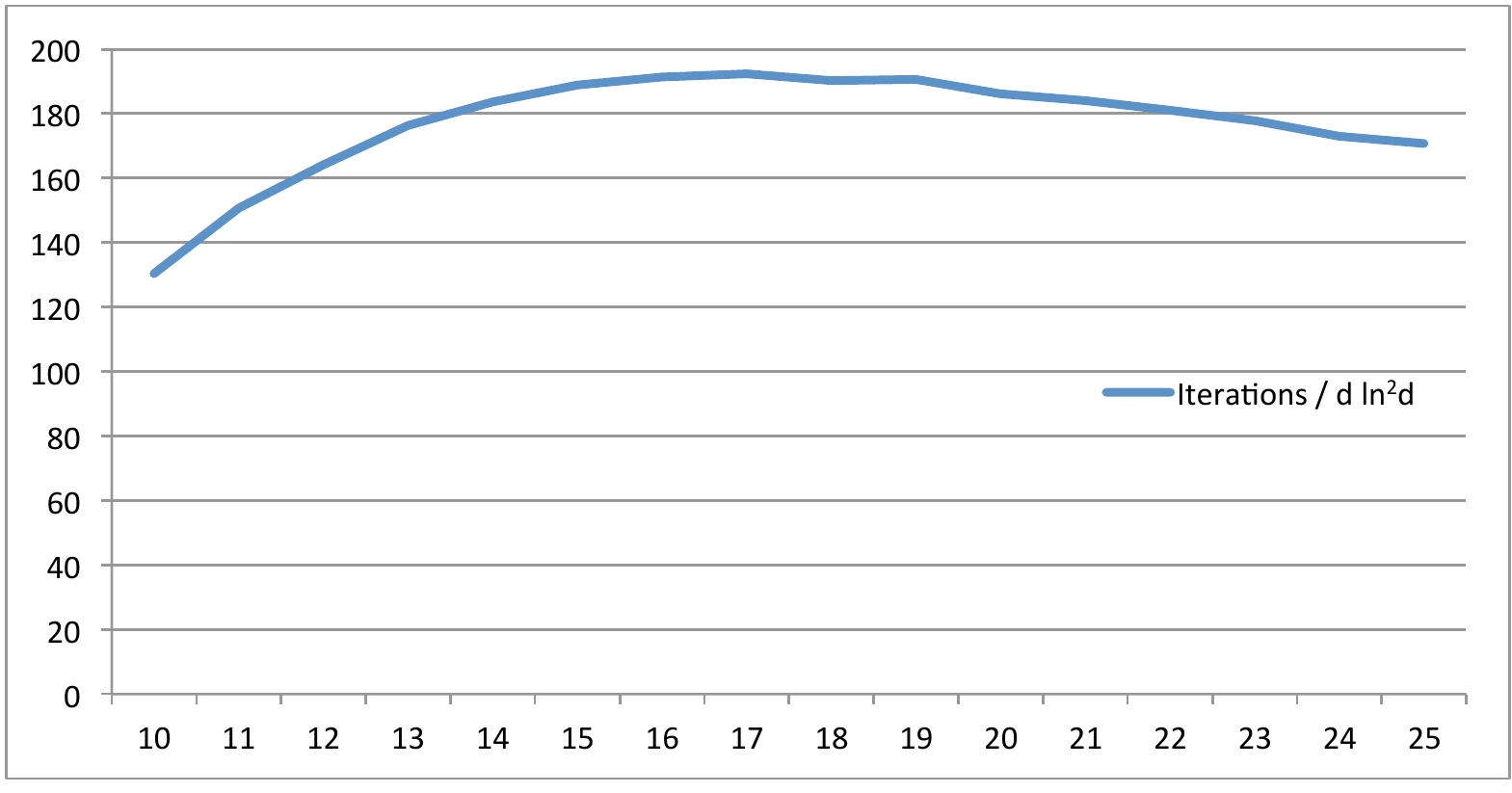}
\caption{Number of Newton iterations for finding centers of hyperbolic components of the Mandelbrot set, divided by $d\ln^2 d$ (column~G in Figure~\ref{Fig:MandelbrotComplexityTable}) vs.\ period. Evidently the number of required iterations scales better than $200\, d\ln^2d$.}
\label{Fig:MandelbrotIterationComplexity}
\end{figure}

\begin{figure}[htbp]
\includegraphics[trim=58 530 94 75,clip,width=.8\textwidth]{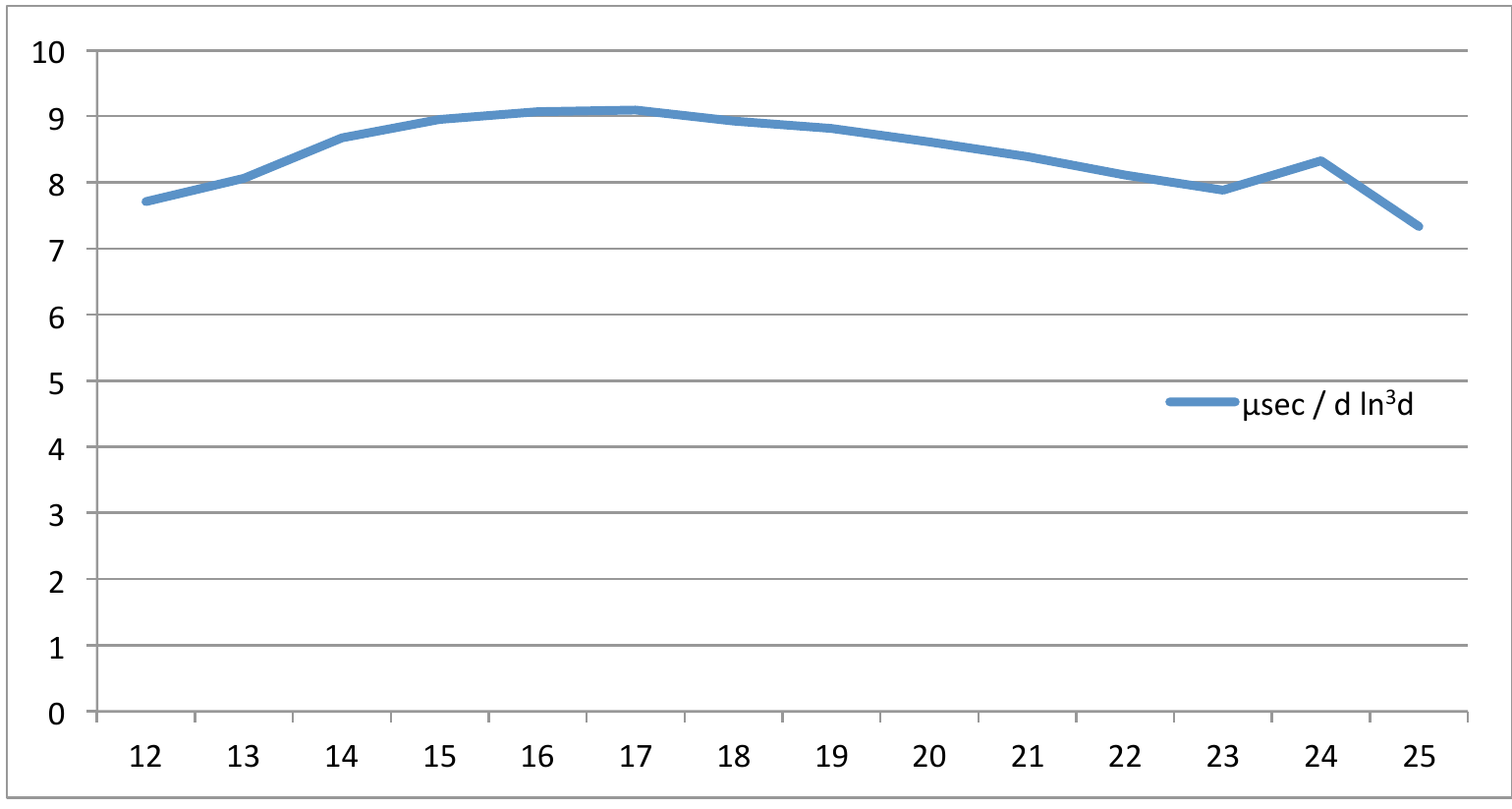}
\caption{Complexity of finding all centers of hyperbolic components of the Mandelbrot set, in terms of computing time, in microseconds divided by $d\ln^3d$ (column I). Evidently, the time complexity scales better than $10\, d\ln^3 d$ (in $\mu\,\text{sec})$. }
\label{Fig:MandelbrotTimeComplexity}
\end{figure}

Roots of $p_n$ are known as \emph{centers of hyperbolic components of the Mandelbrot set} of period $n$ (except those that are already roots of $p_k$ with $k|n$). We found all roots for $n= 25$ (i.e.\ degree $2^{24}$, greater than 16 million) in less than 160 hours (on a standard PC using a single core). 
The detailed results are tabulated in Figure~\ref{Fig:MandelbrotComplexityTable}. In particular, we plot the complexity in terms of required Newton iterations in Figure~\ref{Fig:MandelbrotIterationComplexity}: the diagram clearly shows that the complexity in terms of Newton iterations scales better than $200\, d\ln^2 d$.
The complexity in terms of computing time is shown in Figure~\ref{Fig:MandelbrotTimeComplexity}; it seems to scale better than $10\, d\ln^3 d$. This is of course related to the fact that our degree $d$ polynomials can be evaluated in $\ln d$ operations. We realize that this is an untypical advantage of our polynomials. However, for polynomials given in coefficient form, there are fast methods of parallel evaluation at many points that should compensate for much of the complexity gain \cite{MoenckBorodin}, \cite[Section~8.5]{AHU} (these methods are efficient only when the number of evaluation points is comparable to the degree of the polynomial, which is the case for the Newton algorithm because eventually all roots have to be found by their own Newton orbit).

In order to find all roots, a significantly more sensitive refinement threshold of $R=0.0005$ was used (compared to $R=0.05$ used in the other experiments).

We should mention that centers of hyperbolic components of the Mandelbrot set up to $n=24$ (degrees up to $2^{23}$) were also found by Simon Schmitt in a separate set of experiments by a more sophisticated version of our Iterated Refinement Method: his improvements concern especially the heuristics of the refinement step so that all roots are found with fewer or later refinements and hence significantly fewer applications of the Newton iteration. These are supposed to be published separately.

\begin{figure}[htbp]
\includegraphics[trim=60 530 200 75,clip,width=.55\textwidth]{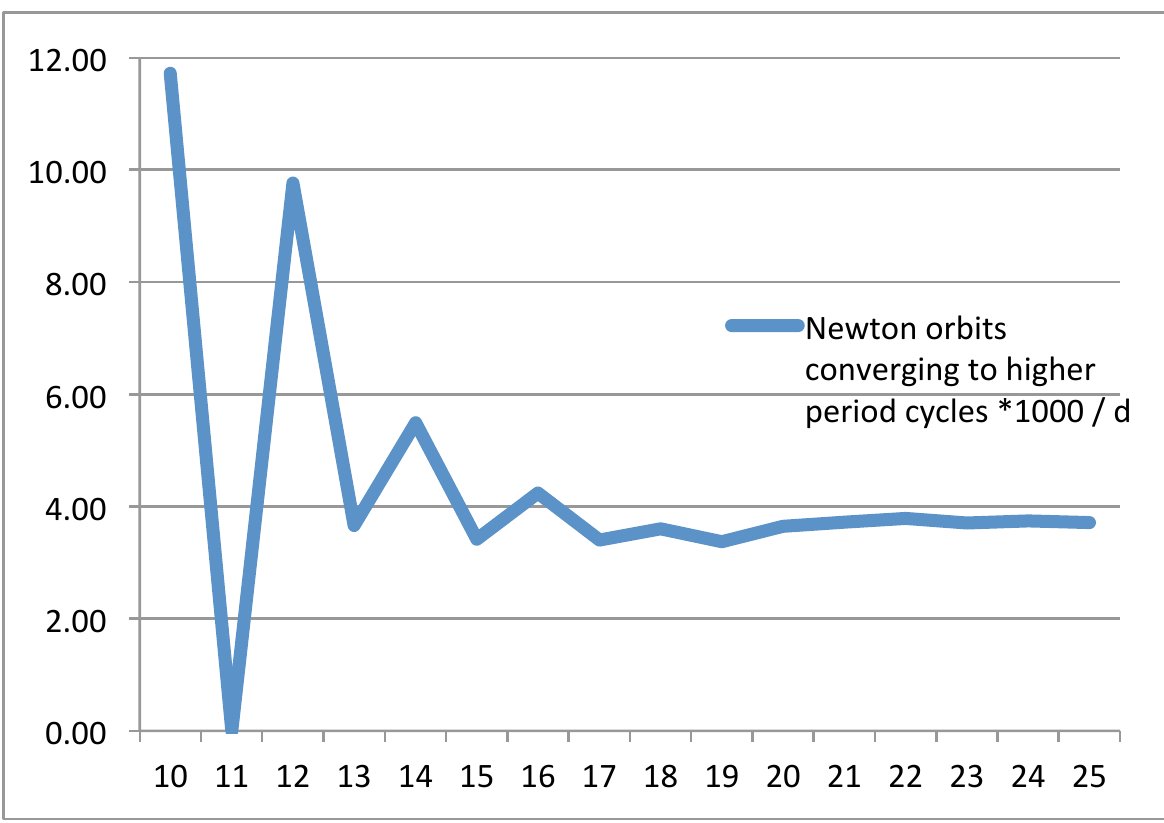}
\caption{Number of Newton orbits converging to attracting cycles of higher period, times $1000/d$: for large $d$, this number seems to stabilize near $0.0037\,d$. }
\label{Fig:MandelbrotHigherCycles}
\end{figure}

\begin{figure}[htbp]
\framebox{
\includegraphics[width=0.93\textwidth]{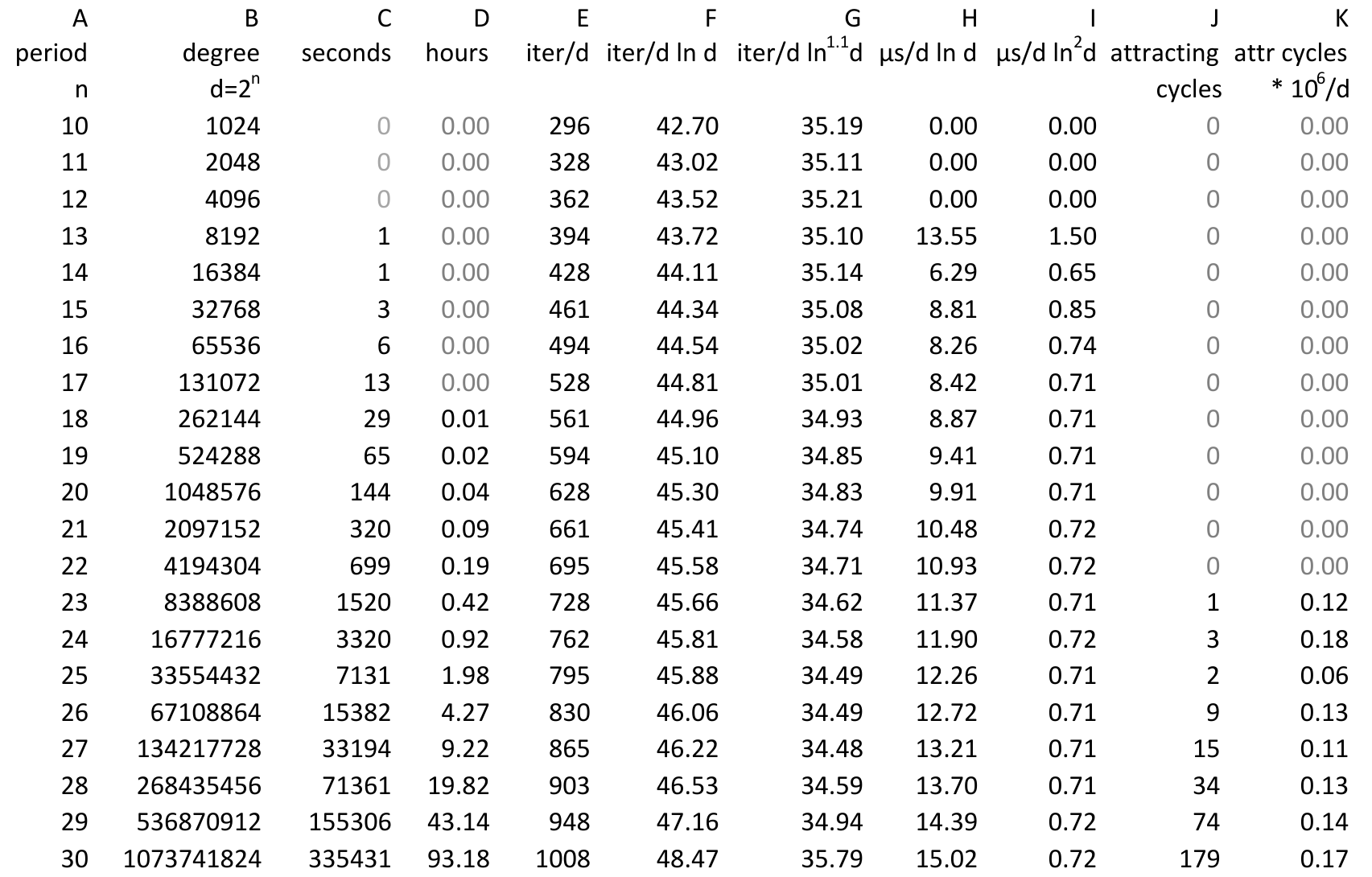}
}
\caption{Finding all periodic points of $p_i(z)=z^2+i$ for period $n\le 30$. The columns are similar as before: period $n$; degree $d=2^n$; computing time in seconds and hours; number of iterations divided by $d$, by $d\ln d$, and by $d\ln^{1.1} d$; computing time (in sec) divided by $d\ln d$, and by $d\ln^{2} d$. The last two columns show the number of orbits that converge to higher period attracting cycles, in absolute numbers and scaled by $d$. (Note that these numbers are \emph{much} smaller than for Mandelbrot centers in Figure~\ref{Fig:MandelbrotComplexityTable}!)
}
\label{Fig:Quadratic_i_ComplexityTable}
\end{figure}

Perhaps not surprisingly, Newton's method occasionally encounters attracting periodic orbits of periods greater than one. The total number of orbits that were detected to converge to attracting cycles of periods $2$ or more is shown in Column J of Figure~\ref{Fig:MandelbrotComplexityTable}; for large $d$, this number seems to stabilize near $0.0037\, d$ (Column K and Figure~\ref{Fig:MandelbrotHigherCycles}). While the total number of such orbits is relatively small, the orbits do not satisfy the easy-to-detect termination condition that a root is found, so the iteration times can be large. Hence, if such orbits are not caught efficiently, they will consume a lot of computing time. It turns out that for the Mandelbrot center polynomials, far more attracting cycles were found than in most other experiments. For comparison, for the other families of polynomials investigated here (periodic points of $z^2+2$ and $z^2+i$), for degree $2^{21}$ no attracting cycles were found at all, and even for degree $2^{30}$ of $z^2+i$ just 179 orbits converged to attracting cycles (see Figures~\ref{Fig:Quadratic2ComplexityTable} and \ref{Fig:Quadratic_i_ComplexityTable}).

\clearpage
\newpage

\subsection{Periodic points of $z^2+2$}

\begin{figure}[htbp]
\framebox{
\includegraphics[width=0.96\textwidth,trim=0 0 0 0,clip]{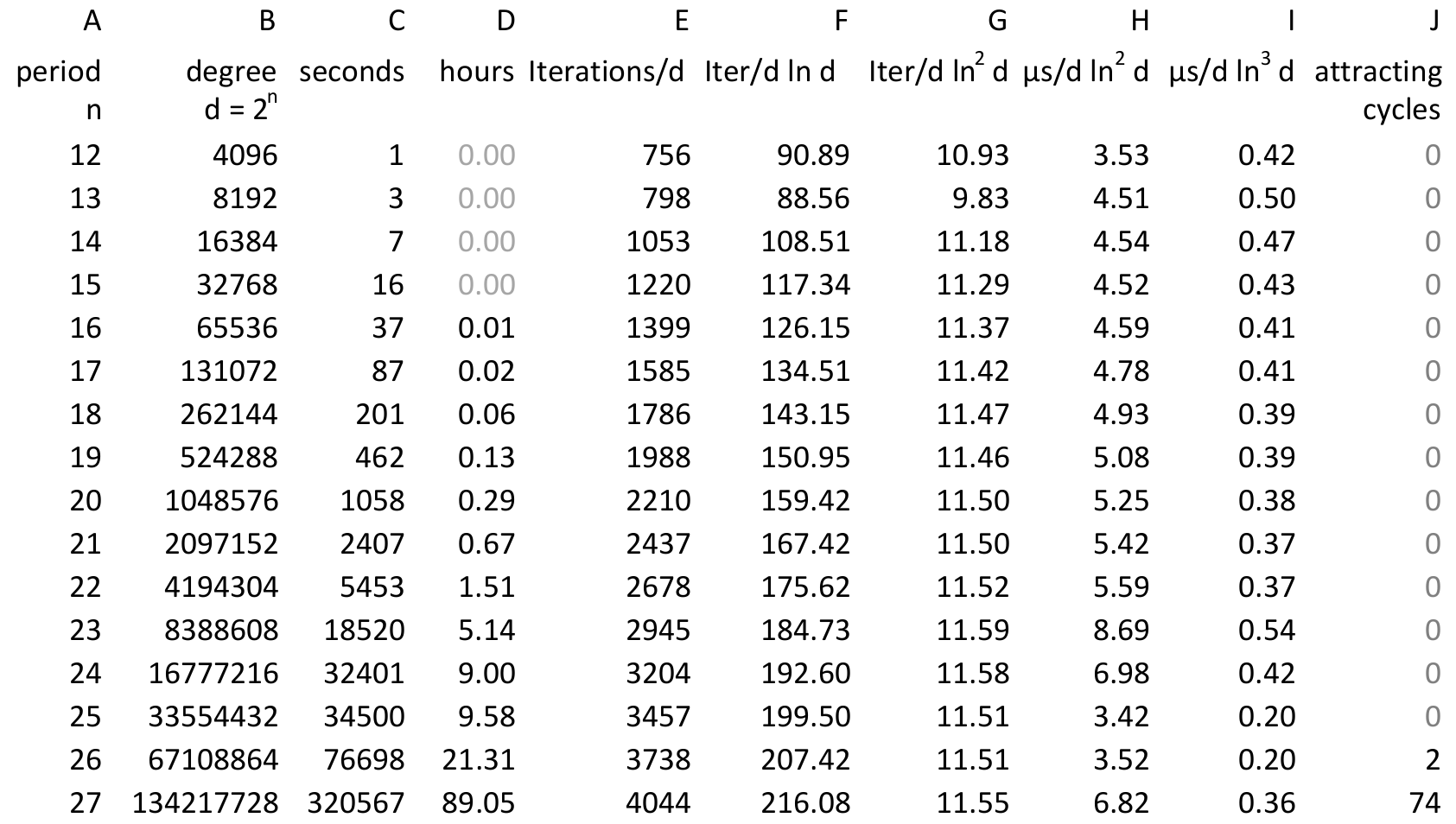}
}
\caption{Finding all periodic points of $p_2(z)=z^2+2$ for periods $n\le 27$. The columns are similar to those in Figure~\ref{Fig:MandelbrotComplexityTable}. The observed complexity in terms of number of iterations is better than $O(d\ln^2 d)$ (see Figure~\ref{Fig:Quadratic2IterationComplexity}), and in terms of computing time approximately $O(d\ln^3 d)$ (with some fluctuations; see Figure~\ref{Fig:Quadratic2TimeComplexity}). The refinement threshold is $R=0.05$, all roots were found. The last column shows the number of orbit that converged to attracting cycles of higher period. }
\label{Fig:Quadratic2ComplexityTable}
\end{figure}

\begin{figure}[htbp]
\includegraphics[width=0.75\textwidth,trim=60 540 170 80,clip]{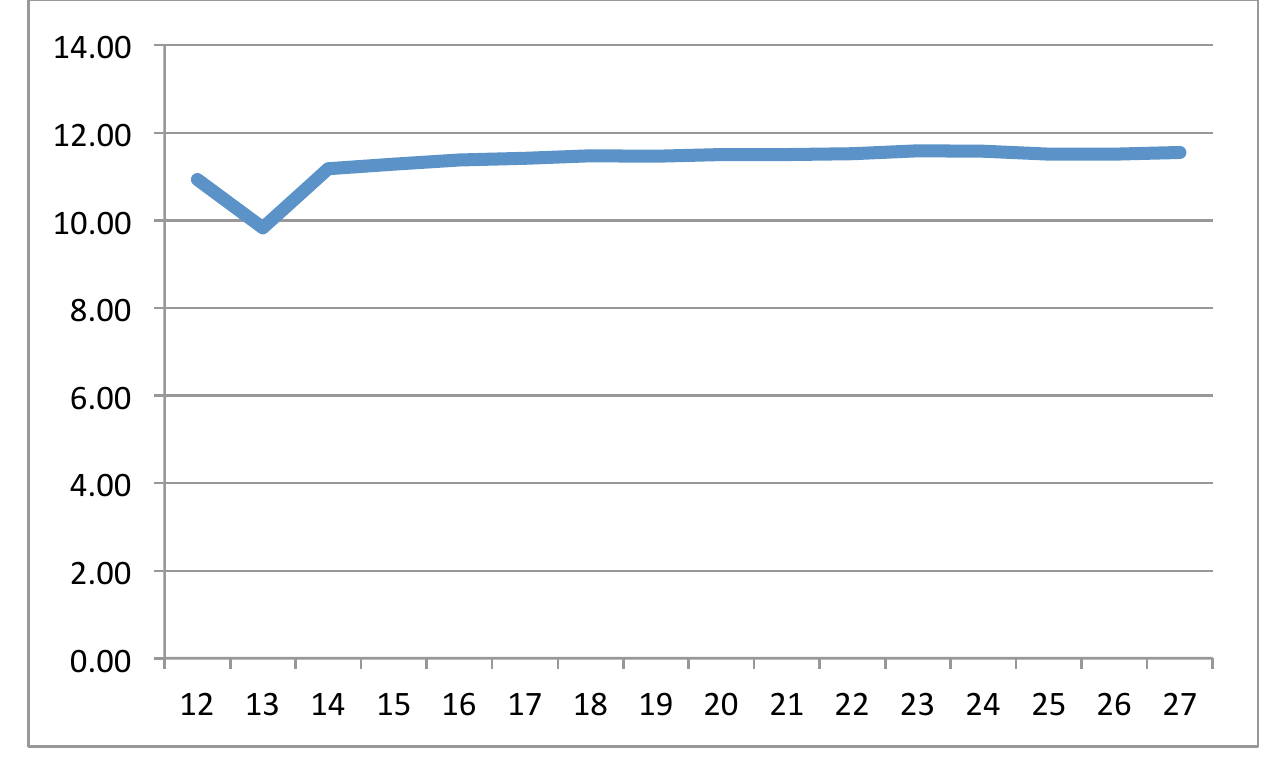}
\caption{Number of iterations required for finding all periodic points of $z^2+2$ for periods $n\le 27$ (divided by $d\ln^2 d$; column G in Figure~\ref{Fig:Quadratic2ComplexityTable}). The complexity in terms of number of iterations seems to scale with $d\ln^2 d$ for periods $n\le 27$, i.e.\ degrees up to $2^{27}$ (greater than $134$ million).  }
\label{Fig:Quadratic2IterationComplexity}
%\end{figure}
%
%\begin{figure}[p]
\includegraphics[width=0.75\textwidth,trim=50 545 155 70,clip]{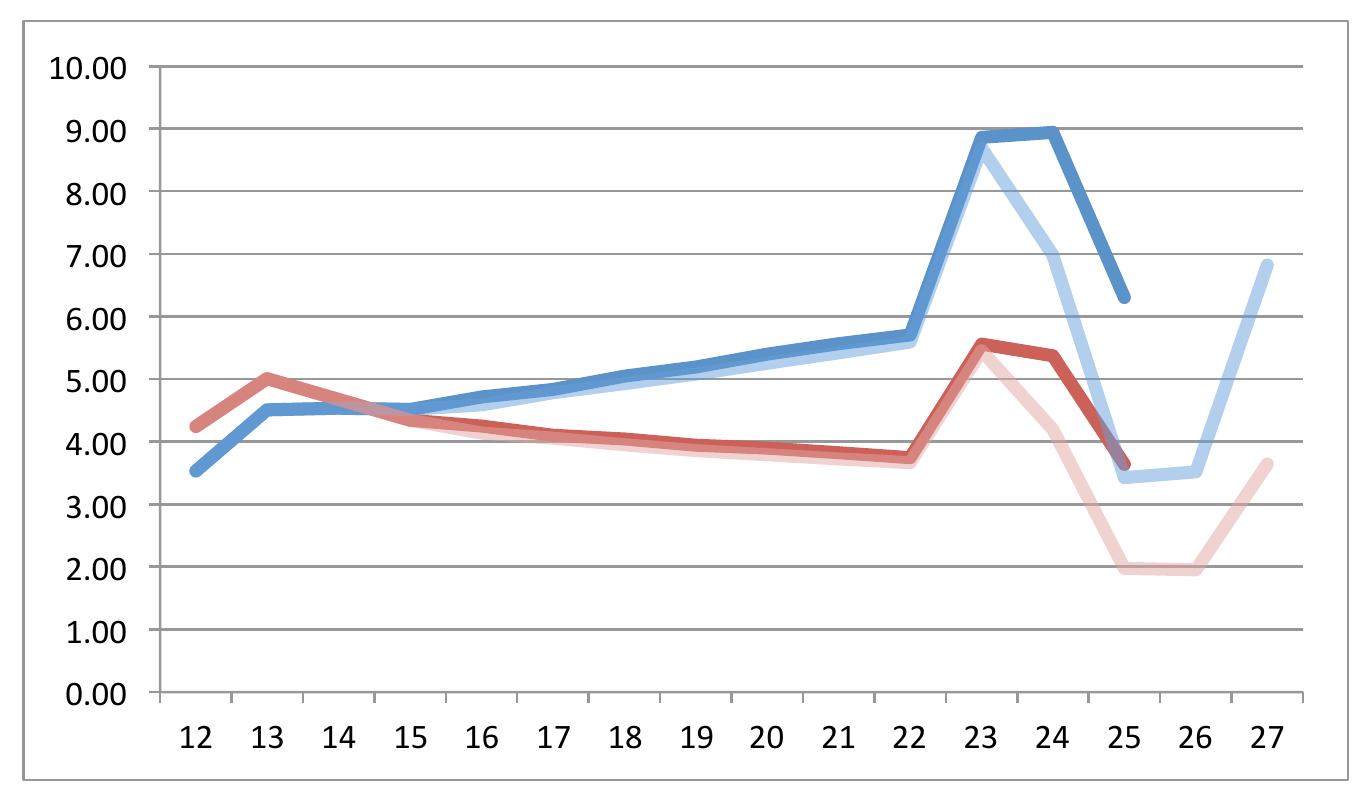}
\caption{Iteration time required for finding all periodic points of $z^2+2$ for periods $n\le 27$, divided by two scale functions. Blue curve: time divided by $d\ln^2 d$ (column H in Figure~\ref{Fig:Quadratic2ComplexityTable}); red curve: time divided by $d\ln^3 d$ (column I), multiplied by $10$ to adjust the scale. The time complexity seems to scale around $0.4\,d\ln^3 d$ (in $\mu\,\text{sec}$, with some fluctuations). Two different runs of the same experiment are shown in different shades of color (the number of iterations remained the same).}
\label{Fig:Quadratic2TimeComplexity}
\end{figure}

For the quadratic polynomial $p_2(z)=z^2+2$, a periodic point of period $n$ is a root of $p_2^{\circ n}(z)-z$, where $p_2^{\circ n}$ stands for the $n$-th iterate of $p_2$. There are exactly $2^n$ periodic points of period $n$ (or dividing $n$). We found all periodic points of period $n\le 27$, i.e.\ for degree up to $134$ million, using a refinement threshold $R=0.05$. The largest degree polynomial took about 89 hours. The results are given in detail in Figure~\ref{Fig:Quadratic2ComplexityTable}. The complexity in terms of required Newton iterations is approximately $2.67 d\ln^2 d$ (see Figure~\ref{Fig:Quadratic2IterationComplexity}). The computing time (shown in Figure~\ref{Fig:Quadratic2TimeComplexity}) is on the order of $0.05\,d\ln^3d $ for degrees $d$ up to 134 million. However, it oscillates significantly by almost a factor of $2$. Since we could not detect a clear reason for this oscillation, we ran the experiment twice. The outcome is very similar, except for large periods. These oscillations have no impact on the overall results of our experiments (but should be investigated in future experiments).
 
Note that this experiment would be outright impossible when expressing $p_2^{\circ n}$ in coefficient form: the constant coefficient alone would have magnitude greater than $2^{2^{n}}$, so for $n=27$ greater than $2^{134\,000\,000}$. Just storing this number would require about 15 megabytes per coefficient (note that good relative precision is not sufficient because all the large terms are subtracted eventually), and we have to accommodate many million coefficients of possibly different magnitudes. 

\clearpage
\newpage

\subsection{Periodic points of $z^2+i$}

We performed the same experiments for $p_i(z)=z^2+i$. The outcome was even more successful in the sense that all roots were found for periods up to $n=30$, that is for degrees up to $2^{30}>10^9$. The computation was very fast: for degree $2^{30}$ it took less than 4 days of computing time on our (not very recent) desktop computers.  It is quite remarkable that a very low computational complexity was observed: only $O(d \ln^{1.1} d)$ for the number of iterations and $O(d\ln^{2})$ for the computing time, for the entire range of degrees up to $2^{30}$ (and notably with small constants in front of the complexity estimates).

Due to the closer mutual proximity of the roots for increasing degrees up to $2^{30}$, the numerical limits of the \texttt{long double} data type seemed to be reached. For that reason, higher numerical precision was used when $|N_p(z)-z|<10^{-13}$, and a root was assumed to be found when $|N_p(z)-z| < \varepsilon_\text{stop} = 10^{-18}$.

\begin{figure}[htbp]
%\framebox
{
\includegraphics[width=0.93\textwidth,trim=0 0 0 0,clip]{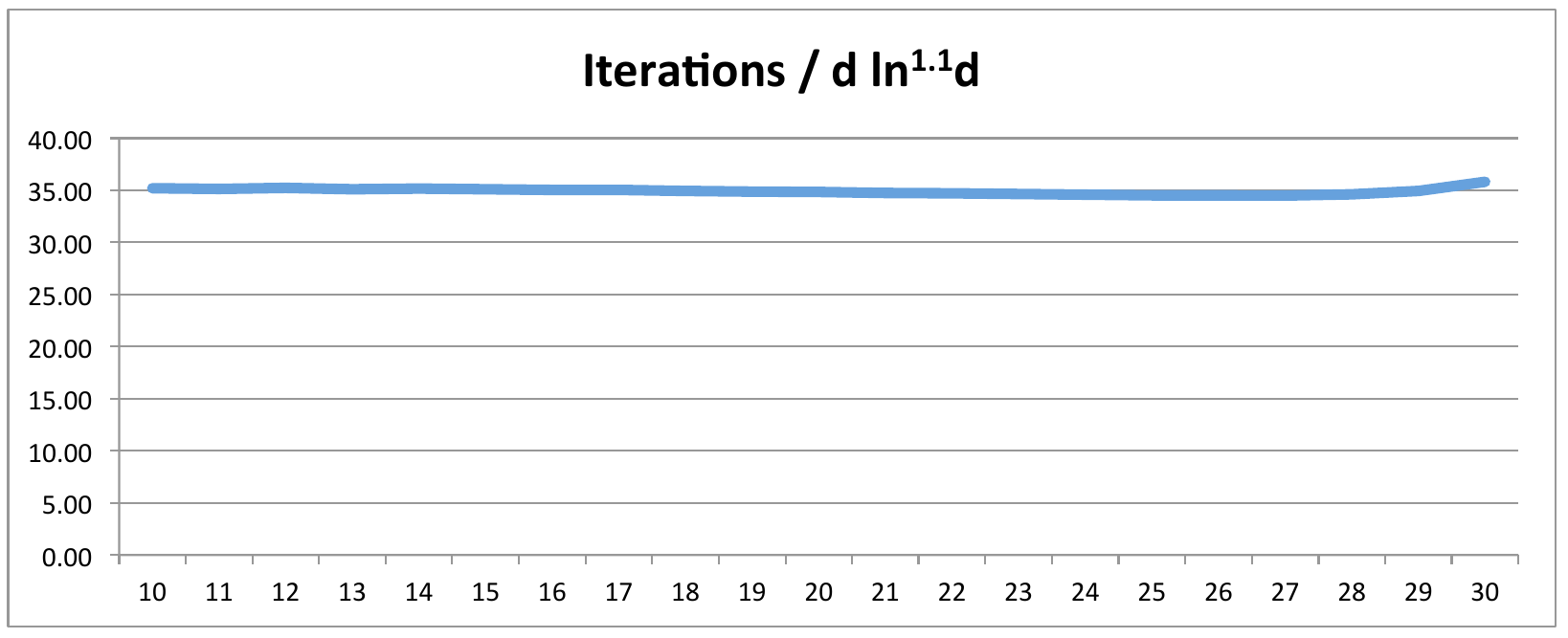}
}
\caption{Number of iterations required to find all periodic points of $z^2+i$ for periods $n\le 30$, divided by $d\ln^{1.1}d$. }
\end{figure}

\begin{figure}[htbp]
{
\includegraphics[width=0.9\textwidth,trim=0 0 0 0,clip]{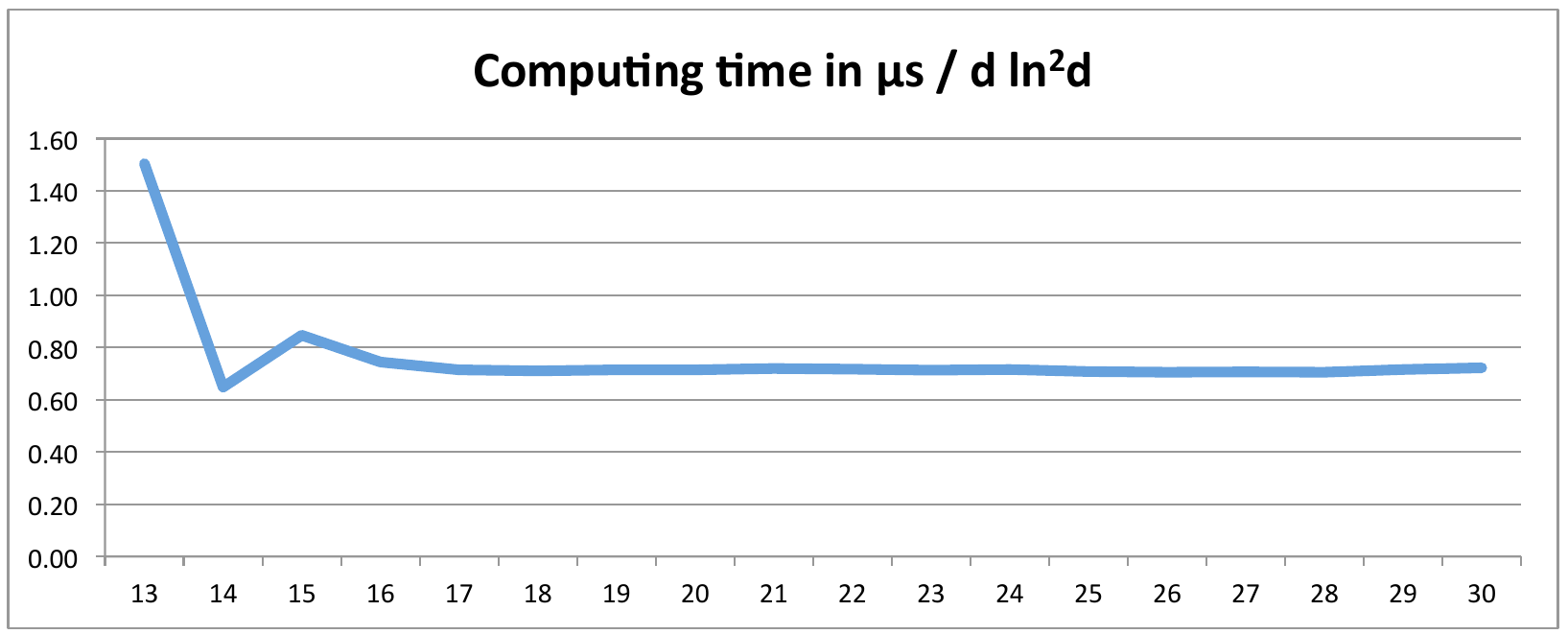}
}
\caption{Computing time necessary to find all periodic points of $z^2+i$ for periods $n\le 30$, divided by $d\ln^{2}d$. }
\end{figure}

\subsection{Computing time per iteration}

The computing time seems to oscillate somewhat unexpectedly in some of our experiments (see for instance Figure~\ref{Fig:Quadratic2TimeComplexity} for $z^2+2$). A natural explanation might come from the fact that the computers had some other tasks running on separate threads (perhaps on separate cores), and while the actual computations would probably require the same amount of computing time, our memory-intense computations might suffer from unexpected memory swaps to the disk. A re-run of the most suspicious experiments showed of course the same number of Newton iterations (this part is deterministic), but also a rather similar oscillation of computing time as before (also shown in Figure~\ref{Fig:Quadratic2TimeComplexity}).

\begin{figure}[htbp]
\framebox{
\includegraphics[width=0.75\textwidth]{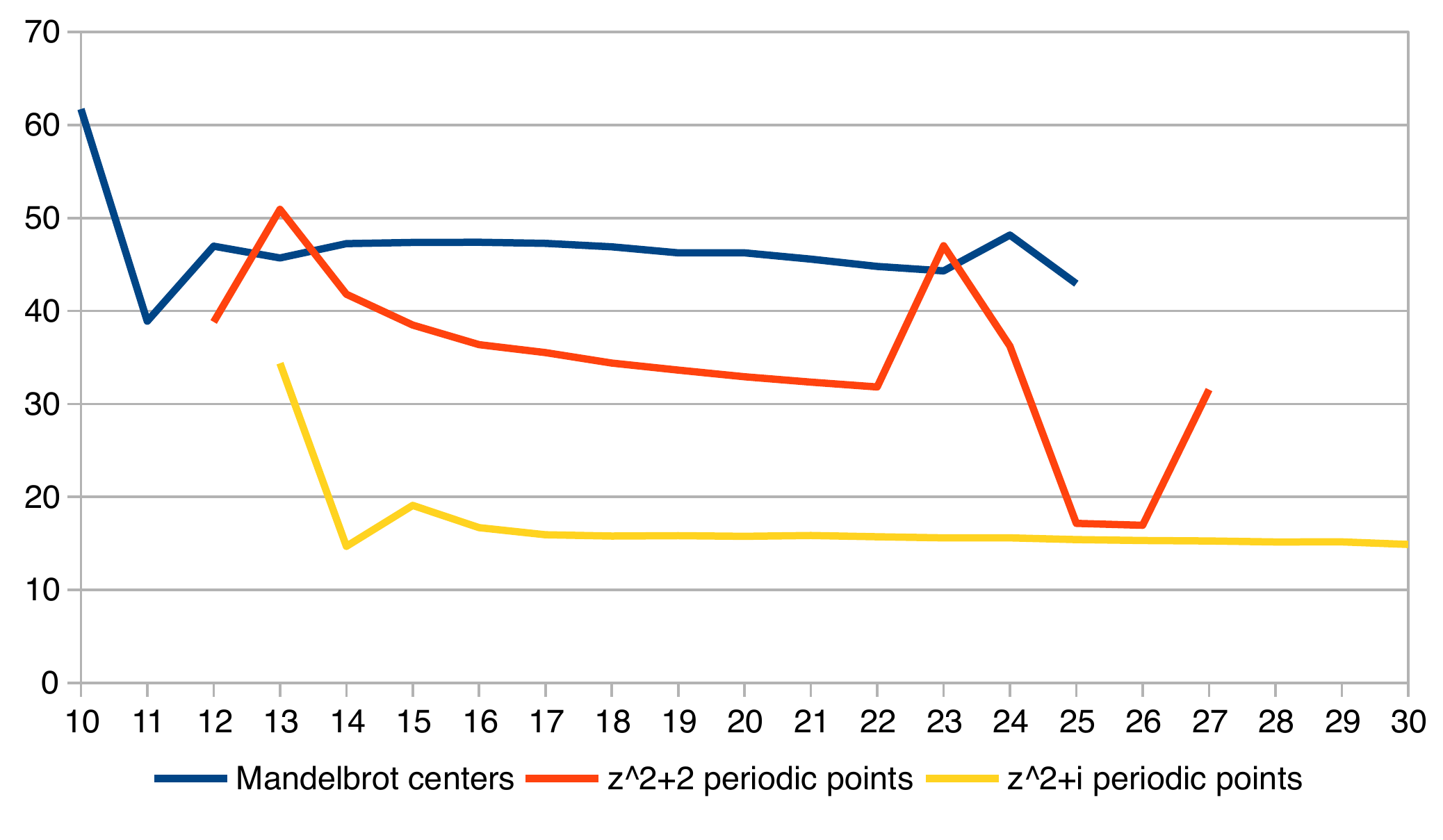}
}
\caption{Computing time in nanoseconds per iteration, divided by $\ln d$ because every evaluation of the polynomial should roughly have complexity $\ln d$; combined for the three sets of experiments. Interestingly, for some experiment these numbers seem to decrease. (This seems to indicate that the Newton iterations do not use most of the computing time, perhaps because the recursion can be implemented particularly efficiently, so other tasks use a definite proportion of the computing time.) Of course, absolute numbers depend on the specific hard- and software used, and not all computations were performed on the same hardware, but the relative behavior should be of structural interest.}
\label{Fig:ComputingTime_rescaled}
\end{figure}

A first step to analyze this behavior consists of comparing computing time per number of Newton iteration steps. Since every Newton iteration for our degree $d$ polynomials requires complexity $O(\ln d)$, we show in Figure~\ref{Fig:ComputingTime_rescaled} computing time divided by $\ln d$. Not surprisingly, this is roughly constant for most series of experiments. Interestingly, for the periodic points of $z^2+i$ one observes a certain decrease in $d$. This might be explained by the fact that not most of the computing time is spent on the Newton iterations, but other tasks require a significant fraction of the time.

Of course, all computing times specified depend on the computer and its software used, in particular efficiency of the compiler and the number of cores available. Therefore, the absolute scales are machine dependent and not very relevant; however, we believe that the relative behavior for changing degrees is of structural interest. The issue of computing time oscillations warrants further study (but stays within a factor of $2$, often less, and is thus not relevant to the key findings of our experiments).

\subsection{What computers were used?}

The computations were run at different locations on different computers from the Intel Xeon family. 
The centers of hyperbolic components of the Mandelbrot sets, as well as periodic points of $z^2+2$, were computed on computers at the university of Bayreuth running at speeds between 2.5 and 3.1 GHz, always using only a single core (of 4 available cores). The periodic points of $z^2+i$ up to period $30$ were computed on computers at Jacobs University running at 2.60 GHz; up to 16 cores were available.

\section{Guaranteeing that all roots are found}
\label{Sec:AllRootsFound}

\subsection{Guaranteeing that all roots are found}
\label{Sub:Guaranteeing}

We admit that our method is heuristic: we do not have an a-priori guarantee that all roots will be found by the Iterated Refinement Newton Algorithm. However, there are several possibilities for verifying a posteriori that all roots have indeed been found. Several of them have been described in detail in \cite[Sections~2.3 and 2.4]{NewtonRobin1}, so we only briefly mention them here; an additional method based on the Newton identities is described in more detail. 

An easy observation shows that for any polynomial $p$ of degree $d$ with associated Newton method $N_p(z)=z-p(z)/p'(z)$, and for arbitrary $z\in\C$, at least one root is contained in the disk $D_s(z)$ with 
\begin{equation}
s=d\cdot|N_p(z)-z|
\;,
\label{Eq:NewtonAccuracy}
\end{equation}
 so small Newton displacement steps $N_p(z)-z$ occur only near the roots. The factor of $d$ can be significantly improved if the approximate locations of many roots are already known \cite[Lemma~4]{NewtonRobin1}. Therefore, if we have $d$ points (coming from $d$ different Newton orbits) that all come with small disjoint disks containing at least one root, then they together describe all roots. This method worked in practice for all degree $2^{20}$ polynomials investigated in \cite{NewtonRobin1}, but it requires that there are no near-multiple roots. For even higher degrees, this guaranteed test becomes increasingly difficult: the factor $d$ in \eqref{Eq:NewtonAccuracy} increases, and the distance between nearby roots decreases.  

An independent and simpler test can be performed based on the Vi\`ete formulas: the coefficients of the polynomial are the elementary symmetric functions of the roots and encode for instance the sum of all roots or their product, and thus provide independent tests of success. More systematically, one can easily compute the value of the power sums $a_k:=\sum_{i}\alpha_i^k$ over all roots $\alpha_i$ directly from the polynomial $p$; this approach is based on the well-known Newton identities that relate the power sums to the coefficients of $p$. Comparing these with the power sums of all roots found is a natural generalization of the Vi\`ete tests. We describe this in detail below (Section~\ref{Sub:NewtonIdentities} for the Newton identities and Section~\ref{Sub:NewtonIdentitiesPractice} for the precision estimates); the overall results are as follows.

In all experiments we describe, we found all roots according to the stopping and distinction criteria described in Section~\ref{Sec:IteratedRefinementNewton}, and then used the power sum criterion to verify the roots found. For all our polynomials, even of maximal degrees, the power sum tests for exponents $x^1$ until $x^{19}$ were run successfully and give reason to conclude that Newton's method with our ad-hoc stopping criterion described above (requiring that $|N_p(z)-z|<\eps_\text{stop}$) found all roots with a typical accuracy of no more than $3\cdot 10^{-16}$.

\subsection{Newton identities}
\label{Sub:NewtonIdentities}
Here we describe a systematic method for checking whe\-ther indeed all roots have been found, based on the well known Newton identities (also known as Newton-Girard formulas) that relate the coefficients of a polynomial (which are up to sign the elementary symmetric functions of its roots) to the power sums of all the roots. This method also makes it possible to locate missing roots in case not all were found (see Section~\ref{Sec:MissingRoots}).

The idea is to use these identities to compute from the coefficients of the polynomial $p$ the sums of the powers of all roots of $p$, then to subtract from these the sums of powers of all roots found: the closer this difference is to $0$, the better the accuracy of the roots found.

Of course it is difficult to find all coefficients of $p$ (and many of them often have forbiddingly large absolute values); but in order to compute the first $m$ power sums of the roots, only the top $m$ coefficients are required (beyond the leading coefficient that we usually scale to $1$ anyway). For large degrees, even the top few coefficients can become very large in practice, and this may present numerical challenges.

We start by writing $p(z)=\prod_{i=1}^d(z-\alpha_j)=\sum_{k=0}^d c_k z^{d-k}$. We have
the leading coefficient $c_0=1$, and the subsequent coefficients $c_k$ are (up to sign) the elementary symmetric polynomials in the roots
\[
c_1 = -\sum_i \alpha_i\;,\quad c_2=\sum_{i< j} \alpha_i\alpha_j\;, \quad c_3=-\!\!\!\sum_{i<j<k}\alpha_i\alpha_j\alpha_k \;, \quad \text{etc.}
\]
We need the power sums 
\[
a_k:=\sum_{i=1}^d\alpha_i^k
\] 
for $k\le m$. They can be computed from the coefficients $c_k$ by the Newton identities as follows:
\begin{align*}
-a_1 &= c_1 \\
-a_2 &= c_1a_1+2c_2 \\
-a_3 &= c_1a_2+c_2a_1+3c_3 \\
-a_4 &= c_1a_3+c_2a_2+c_3a_1+4c_4
\end{align*}
\hide{
\begin{align*}
-c_1&= a_1 \\
-2c_2&= c_1a_1+a_2 \\
-3c_3&=  c_2a_1+c_1a_2+a_3 \\
-4c_4 &= c_3a_1+c_2a_2+c_1a_3+a_4
\end{align*}
}
and so on for subsequent coefficients. Therefore, each $a_k$ can be computed recursively from the $c_k$ and the previously computed $a_i$.

In our applications, the polynomial $p$ is not given in coefficient form, so in the first step the coefficients $c_k$ have to be computed. Only the coefficients $c_1\dots c_m$ are required. They can be found easily using the recursion that defines $p$: if in the recursion in every step only the $m$ top-most coefficients are kept, then the final result will correctly yield the desired coefficients $c_1,\dots, c_m$. 

To compute the coefficients $c_k$, one needs $\ln d$ recursion steps, each with complexity $m^2$ (or less if more efficient multiplication is implemented), for a total of $O(m^2\ln d)$. The theoretical power sums $a_k$ for $k\le m$ are computed through the Newton identities, which are a triangular system of linear equations of dimension $m$, so they are evaluated in complexity $O(m^2)$. Both steps are negligible compared to the computation of the actual power sums from the roots found, which has complexity $O(dm)$.

%\clearpage\newpage

\subsection{Newton identities in practice and precision estimates}
\label{Sub:NewtonIdentitiesPractice}

\begin{figure}[htbp]
\includegraphics[width=0.8\textwidth,trim=62 530 67 75,clip]{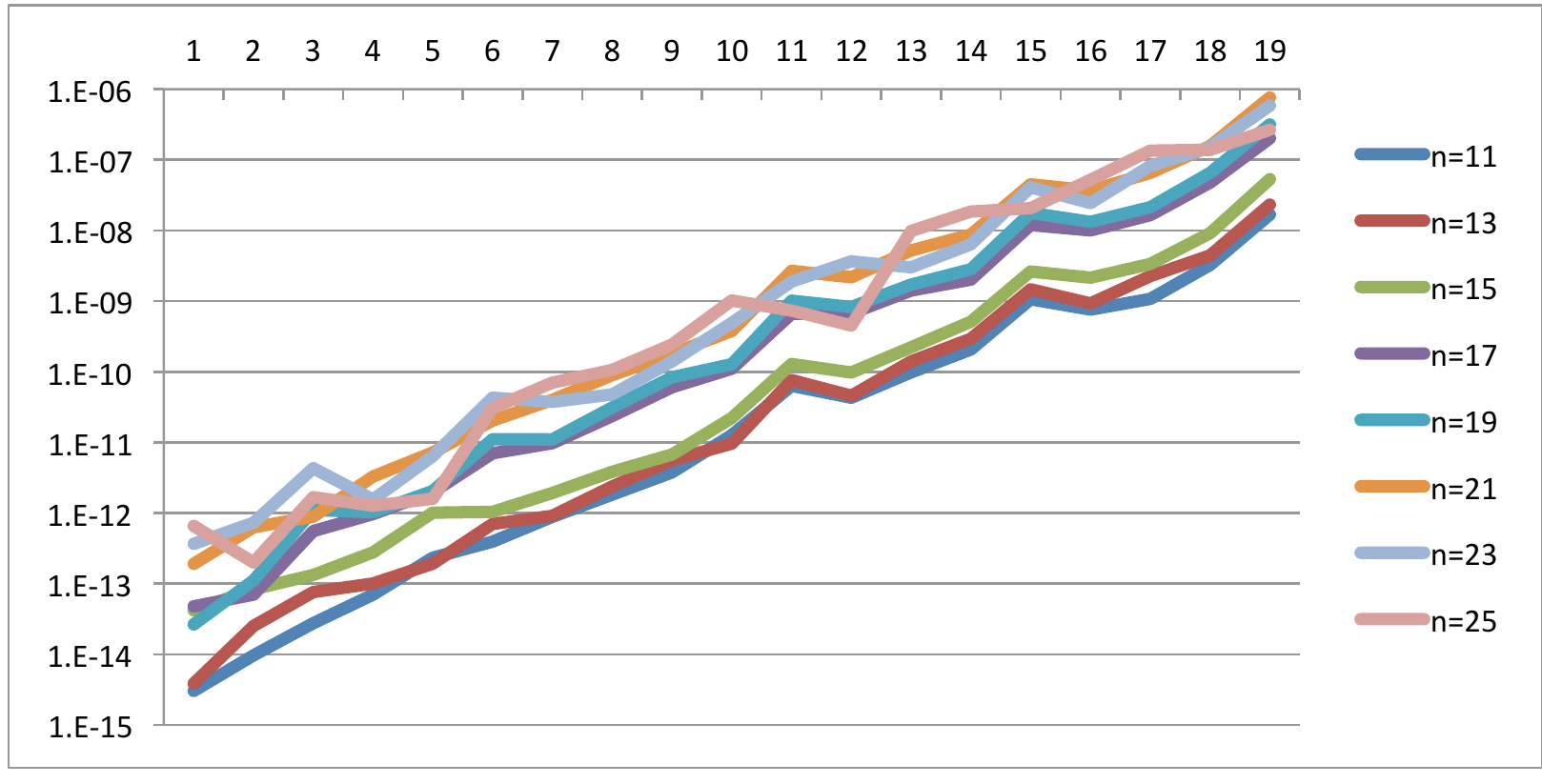}
\begin{picture}(0,0)
\put(-250,107){Centers of the Mandelbrot set}
\end{picture}
\includegraphics[width=0.8\textwidth,trim=65 575 110 75,clip]{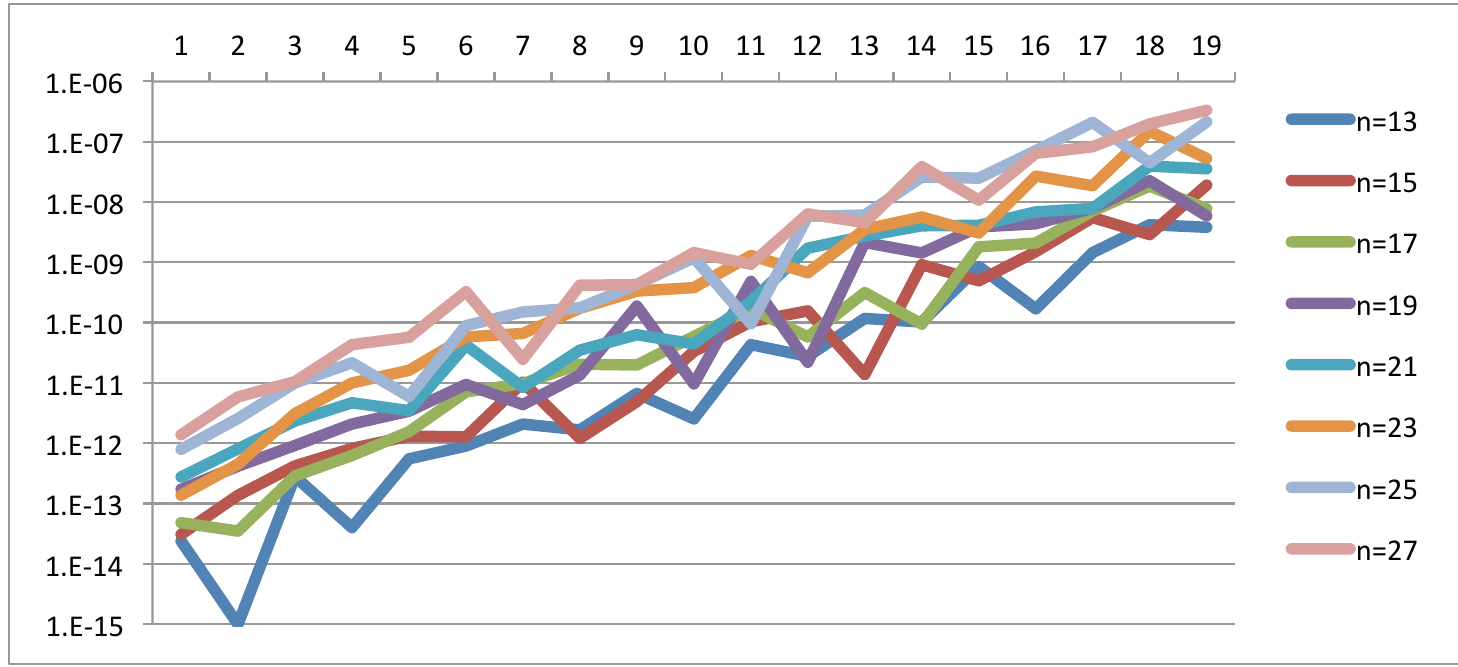}
\begin{picture}(0,0)
\put(-250,94.5){Periodic points of $z^2+2$}
\end{picture}
\includegraphics[width=0.81\textwidth,trim=-3 0 1 0]{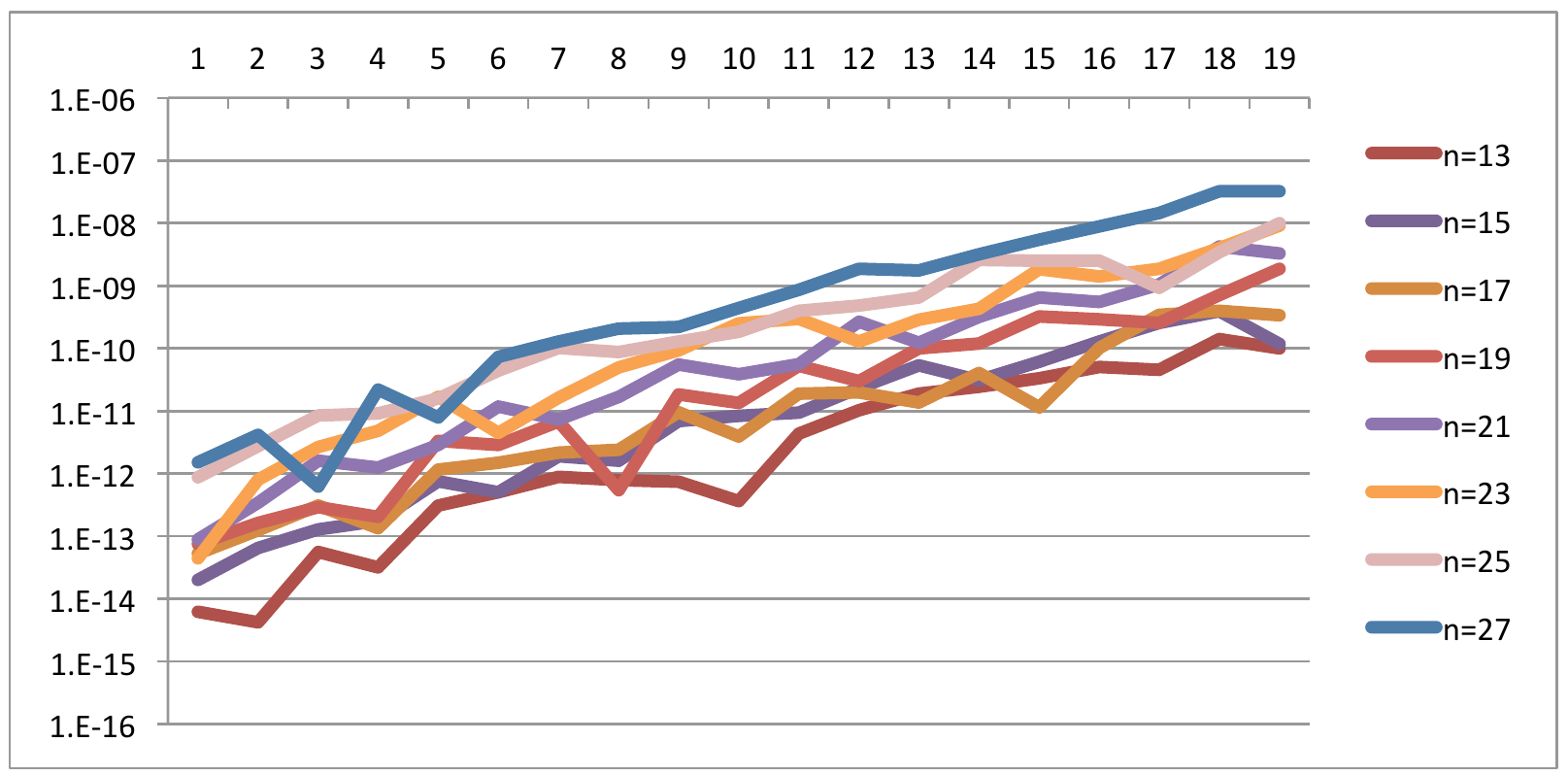}
\begin{picture}(0,0)
\put(-250,104){Periodic points of $z^2+i$}
\end{picture}
\caption{A posteriori verification that all roots were found, using the Newton identities for exponents $k\le 19$ (horizontal scale), for odd periods $n$. The vertical value is the deviation of the power sums of all roots from the (Gaussian) integer value predicted by the Newton identities. Top: centers of the Mandelbrot set; middle: periodic points of $z^2+2$; bottom: periodic points of $z^2+i$.  \hide{\reminder{Redo the graph for $z^2+i$ for $n\le 30$ when data are there; adapt caption.}}
}
\label{Fig:NewtonIdentityPrecision}
\end{figure}
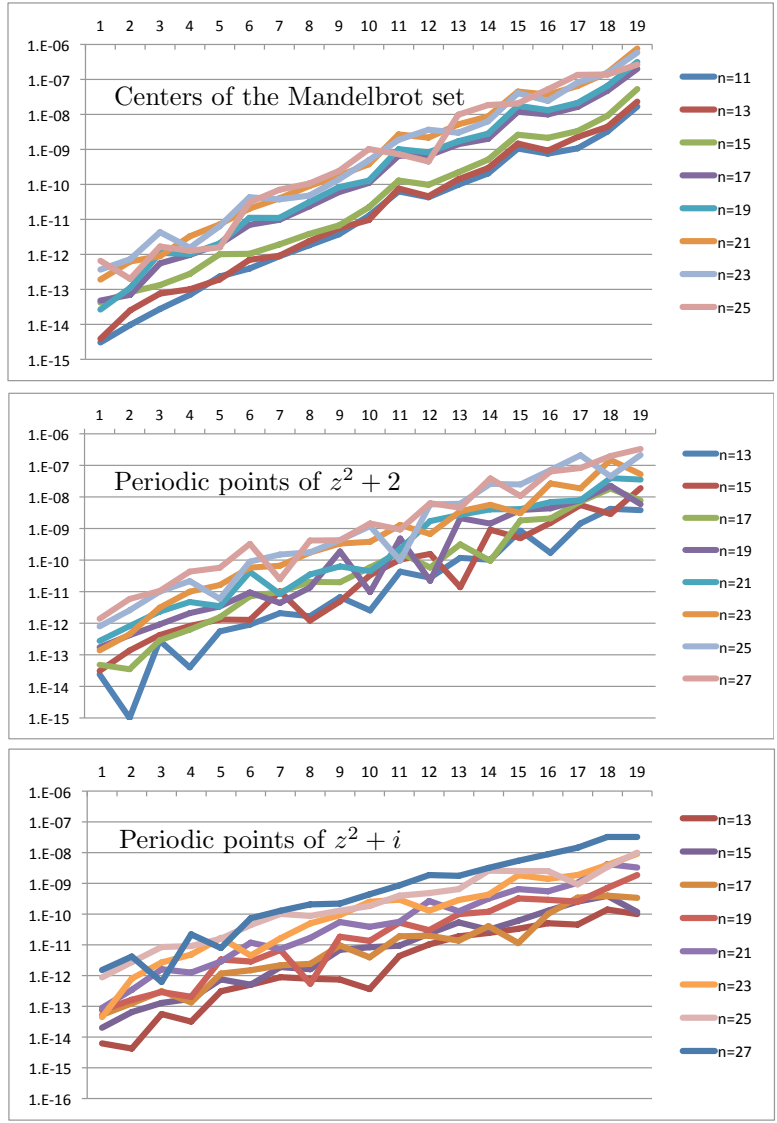

For the three families of polynomials, we computed the power sums for exponents $k\le 19$, and compared them to the actual values computed using the Newton identities. The differences are shown in Figure~\ref{Fig:NewtonIdentityPrecision}. 

The observed errors are between $10^{-15}$ and $10^{-8}$, not surprisingly with increasing errors for higher degrees and higher powers. Within reasonably expected numerical error (discussed below), we see these computations as convincing confirmation that indeed all roots were found. 

One way to model the error in these tests is to assume that all roots found have numerical errors of average size $\delta$, distributed equally and independently. Specifically for the power sums with exponent $k=1$, the sum is then a random walk in $\C$ of $d$ steps with average step size $\delta$, so the accumulated error should be of size $\delta\sqrt d$. 
As an example, for the $d=2^{27}$ periodic points of $z\mapsto z^2+2$ of period $27$, the numerically computed power sum for $x^1$ was $-1.388\cdot 10^{-12}+2.220\cdot 10^{-15}i$
with a deviation of less than $1.4 \cdot 10^{-12}$ from the true integer value $0$ computed via the Newton identities. Our random walk assumption implies $\delta\sqrt d\approx 1.4 \cdot 10^{-12}$, so we can estimate $\delta\approx 1.4\cdot 10^{-12}\cdot 2^{-13.5}\approx 1.2\cdot 10^{-16}$.

\begin{figure}[htbp]
\includegraphics[width=0.8\textwidth,trim=55 490 150 75,clip]{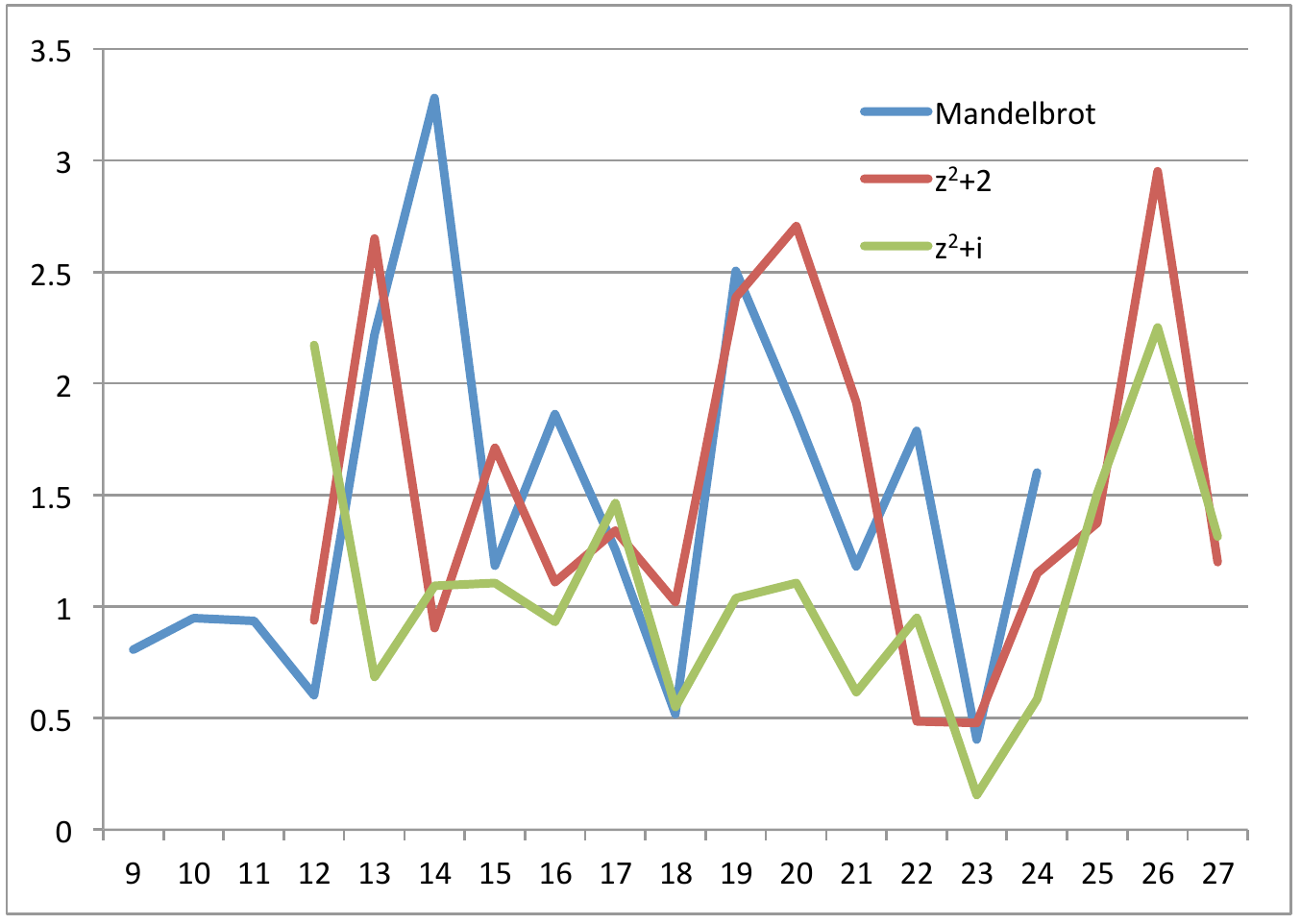}
\caption{Expected precision of all roots, in units of $10^{-16}$, for the three series of experiments (vertical). The expected error is computed from the $x^1$ terms of the Newton identities, based on the assumption of independent and equal distribution (Brownian motion in $\C$). The horizontal scale shows $N$ so that the degree is $2^N$ (the periodic points have period $N$, and the Mandelbrot centers have period $N+1$). Observe that all values are of approximately constant size, mostly within a factor of~$5$ from each other. These data strongly suggest that all roots have been found with typical accuracy of $\delta\approx 3\cdot 10^{-16}$ or better. \hide{\reminder{Add data for $n\le 30$!}}
}
\label{Fig:ExpectedError}
\end{figure}

Similar estimates were made for all series of experiments and all power sums $x^1\dots x^{19}$. 
Figure~\ref{Fig:ExpectedError} shows that the results for $k=1$ are consistent with the interpretation that typically $\delta<3\cdot 10^{-16}$.

For $k>1$, we encounter a precision issue for large numbers. Let $\rho$ be the largest absolute values of the roots of a given polynomial; we have $\rho\approx 2$ for the Mandelbrot centers, $\rho\approx 1.74$ for periodic points of $z^2+2$, and $\rho\approx 1.48$ for periodic points of $z^2+i$. Then $|\alpha_i^k|\approx \rho^k$ for some $i$. These numbers increase exponentially with large $k$ and are eventually subtracted from each other, or from the exact coefficients that in our cases are (Gaussian) integers, so (for number formats with fixed lengths mantissae) the absolute errors in the sums increase with $\rho^k$; this would lead to approximately linear graphs in the accuracies shown in Figure~\ref{Fig:NewtonIdentityPrecision}, with slope $\ln\rho$. Indeed, the average slopes are close to $2.26$ (for the centers of the Mandelbrot set), to $1.92$ (for periodic points of $z^2+2$), and to $1.65$ (for periodic points of $z^2+i$); so in all three cases they are comparable to, and somewhat larger, than the predicted lower bounds of the accuracy. In any case, the most accurate tests for the achieved accuracy are those with low exponents $k$.

It might be interesting to report a few new challenges for these computations, coming from finite accuracy. One can use the Newton identities in both directions: from the numerically computed power sums one can compute the coefficients and compare these with the true values from the polynomials, or do the computation backwards starting with the coefficients and compare the power sums. The latter approach, from coefficients to power sums, seems more natural in our case because we do most computations with (Gaussian) integers, so without losing precisions. However, this involves \emph{much} larger numbers (our polynomials have very large coefficients, and we had to use integer arithmetic that can handle products of these without losing valid digits); the computation in the other direction does not have the advantage of integers, but the numbers required are significantly smaller. 
 
Another observation is that the very summation of the root powers loses valid digits. We had mentioned above that especially the high powers of the roots have limited absolute precision. In addition, when adding these large numbers, we often encounter very large sums (consistent with the large predicted values from the Newton identities). For instance, the powers $k=18$ for period $n=27$ of periodic points of $z^2+2$ sum up to about $-2\cdot 10^{11}$, which reduces the available absolute accuracy further. Even when the power sum has small absolute value (such as for $k=19$ for the same polynomial), the intermediate sums might become large before decreasing eventually (especially when the roots found are sorted for instance by increasing real parts). We thus ran the addition a second time modulo $1$ and modulo $i$, hence keeping track only of the fractional parts of the sums (having checked that the integer part of the solution is consistent), so that all available accuracy can measure the difference to the exact integer value. This yielded results consistent with the initial computations, but with substantially greater precision. These improved results were used for the precision checks and are shown in Figures~\ref{Fig:NewtonIdentityPrecision} and \ref{Fig:ExpectedError}.

\section{Finding missing roots}
\label{Sec:MissingRoots}

A structural disadvantage of our iterated refinement Newton root-finding method is that it does not come with an a priori guarantee that it will find all roots, and if in the end it turns out that some are missing, it is not clear where a refinement of starting points will recover the missing roots. This is different from the simpler linear scheme described in \cite{NewtonRobin1}: there we usually start with $4d$ equidistributed points on a single circle, and if some roots are missing, then this initial set of points can get refined for additional orbits. The crucial difference is that in the linear scheme from \cite{NewtonRobin1}, every orbit starts at the same large circle surrounding all roots, so every orbit needs to run $\Omega(d)$ iterations before even getting close to the roots. In our Iterated Refinement Newton Algorithm, most orbits start very close to the roots, so far fewer iterations are required, but the refinement takes place away from the locus where we have good control. We gain a lot of speed at the expense of giving up guaranteed convergence.

Indeed, in some experiments it does happen that the Iterated Refinement Newton Algorithm finds almost all roots, but a few of the more than one million roots are missing (in one case that we describe in detail, three out of $134$ million roots were missing). Roots may be missing either  because the refinement was not sensitive enough, or because the roots were found but not declared sufficiently different from nearby roots that were reported.  
Our global scheme does not have a way to pinpoint where the missing roots were lost, and thus where additional starting points are required to discover these.

Explicit deflation (dividing the original polynomial by the product of all linear factors corresponding to the roots found) is tempting but not an option: it is numerically unstable unless the roots are found in a particular order, and it destroys any structure inherent in the polynomial (and it might have complexity in itself comparable to the entire root finding process).

Here we describe three methods, presumably well known, that have low complexity and that find missing roots efficiently when the number of remaining roots is small: if the given polynomial $p(z)=\prod_{i=1}^d(z-\alpha_i)$ has degree $d$ and $m$ roots are missing, then the complexity of this post-processing is essentially $O(dm)$, up to logarithmic factors. 
We are grateful to Dario Bini, Victor Pan and Michael Stoll for helpful suggestions on these topics.

\subsection{Implicit deflation}
The first method is ``implicit deflation''.
Suppose the roots found are $\alpha_1,\dots,\alpha_{d-m}$ and the missing roots are $\alpha_{d-m+1}$, \dots, $\alpha_d$. Write $q(z)=\prod_{i=d-m+1}^d(z-\alpha_i)$. Of course we do not know $q$ directly.

The idea of implicit deflation is to evaluate Newton's method for $q(z)$ in terms of $p$ and the roots found: we have $q(z)=p(z)/\prod_{i=1}^{d-m}(z-\alpha_i)$ and thus 
\begin{equation}
\frac{q'(z)}{q(z)} = \frac{p'(z)}{p(z)} - \sum_{i=1}^{d-m} \frac{1}{z-\alpha_i} 
\;
\label{Eq:LogDerivativeAdditive}
\end{equation}
(for products of polynomials, the logarithmic derivative is additive).
We can thus evaluate $N_q(z)=z-q(z)/q'(z)$ using only $p'(z)/p(z)$ (which is computed by the usual Newton method for $N_p$) and the roots found so far. 

If $m\ll d$, then the final sum dominates the complexity: every iteration step of $N_q$ has complexity essentially $O(d)$, while the number of necessary Newton iterations depends only on the degree $m$ of $q$. Depending on the efficiency of root finding for this (relatively) low degree polynomial, the number of required Newton iterations is between $O(m\ln m)$ and $O(m^2)$, so the total complexity of finding the remaining roots is between $O(md\ln m)$ and $O(m^2d)$. 

\subsection{Ehrlich-Aberth iteration}
A second method to locate the missing roots is the Ehrlich-Aberth method (which is underlying the successful implementation \texttt{MPSolve} to find all roots).  This method is a parallel iteration in $d$ variables, where the $j$-th component is the Newton method of $p(z)/\prod_{i \neq j}(z-\tilde\alpha_i)$ and the $\tilde\alpha_i$ are the current approximations to all $d$ roots. As with implicit differentiation, one does not need to compute $\prod_{i\neq j}(z-\tilde \alpha_i)$ explicitly for the Newton displacement, but uses an analog of \eqref{Eq:LogDerivativeAdditive}.

All $d$ coordinates are updated simultaneously. If $d-m$ roots are already known, then this method is an iteration in only $m$ variables. In particular, if $m=1$, then implicit deflation equals the Ehrlich-Aberth method with $d-1$ coordinates fixed; these two methods are thus closely related.

\subsection{Newton identities} 
\label{Sec:MissingRootsNewtonIdentities}

A third method for locating missing roots is based on the Newton identities described in Section~\ref{Sub:NewtonIdentities} that relate the coefficients of a polynomial (which are up to sign the elementary symmetric functions of its roots) to the power sums of all the roots. 
This method may be in less frequent use, even though it is based on very classical ideas. It allows one to construct a polynomial $q$ of degree $m$ 
the roots of which are exactly the missing roots of $p$. This polynomial will be constructed in terms of its coefficients using sums of different powers of the roots as well as Newton's identities.

The idea is to use these identities to compute from the coefficients of the polynomial $p$ the sums of the powers of all roots of $p$, then to subtract from these the sums of powers of all roots found: we thus obtain the sums of powers of the missing roots, and from these one can compute the coefficients of $q$ by applying the Newton identities backwards (where $q$ as above is the polynomial that has exactly the missing roots). We describe this method in somewhat more detail.

We described in Section~\ref{Sub:NewtonIdentities} linear relations that allow one to compute from the coefficients of the polynomial $p$ the sums $a_k$ of the $k$-th powers of all its roots. 

Now suppose the $d-m$ roots $\alpha_1\dots,\alpha_{d-m}$ of $p$ have already been computed, and $m$ roots $\alpha_{d-m+1},\dots,\alpha_d$ are missing. We can compute their power sums
\[
b_k=\sum_{i=d-m+1}^d  \!\!\! \alpha_i^k = \sum_{i=1}^d \alpha_i^k -\sum_{i=1}^{d-m} \alpha_i^k = a_k - \sum_{i=1}^{d-m} \alpha_i^k
\]
from the $a_k$ and the $d-m$ roots found. Finally, we are looking for the polynomial 
\[
q(z)=\prod_{i=d-m+1}^d \!\!\! (z-\alpha_i)=\sum_{k=0}^m e_kz^{m-k}
\;,
\]
where the $e_k$ are the coefficients of $q$. The coefficients $e_k$ can be determined from the $b_j$ using the Newton identities backwards: 
\begin{align*}
-\,\,\,\,\!e_1& = b_1 \\
-2e_2&= e_1b_1+b_2 \\
-3e_3 &= e_2b_1+e_1b_2+b_3 \\
-4e_4 &= e_3b_1+e_2b_2+e_1b_3+b_4 \\
& \,\,\; \vdots
\end{align*}
Therefore, once we know $b_1,\dots, b_m$, we can find $q$ in coefficient form and then apply any root finder to this polynomial of degree $m$. 

In order to find the $m$ coefficients $e_1,\dots,e_m$ of $q$, only the power sums $b_1,\dots,b_m$ are required in the Newton identities backwards, and hence only the power sums $a_1,\dots,a_m$ and the coefficients $c_1\dots,c_m$ are required. We described earlier in Section~\ref{Sub:NewtonIdentities} that these coefficients can be computed with little effort even when our polynomials $p$ are not given in coefficient form. 

The complexity of this method has been described in Section~\ref{Sub:NewtonIdentities}: it is dominated by the computation of the power sums $b_1\dots b_m$, which requires $O(dm)$ operations. The remaining computations for $m$ missing roots can be performed in $O(m^2\ln d)$ operations.

One structural advantage of this method is that the original roots are required only for the simple computation of the power sums, not for the subsequent root finding iteration. In our applications, the roots were found on a different computer system than the one on which the missing roots were located, so not all $d-m$ previously found roots had to be transmitted, but only $m$ easily computed power sums (it makes a substantial difference to transmit, or even to store, $2^{27}-3$ complex roots or $3$ power sums). This advantage comes with a certain disadvantage, though: the inaccuracies of all $d-m$ roots contribute equally to inaccuracies in the power sums, no matter how far they are from the missing roots; in practice this limits the achievable accuracy of this method (for the methods of implicit deflation or Ehrlich--Aberth, the inaccuracy of roots that are far from missing roots becomes far less important). 

The two purposes of Newton identity, to find missing roots and to check the roots found and estimate the accuracy achieved, can be combined: after $m$ missing roots have been reconstructed using the first $m$ coefficients, one can compute $m'$ additional coefficients and hence $m'$ additional power sums, and use these to verify that now all roots have been found with great accuracy. This will be illustrated in an example in Section~\ref{Sub:NewtonIdentitiesPractics}.

%\newpage

\subsection{Newton identities in practice: locating three missing roots}
\label{Sub:NewtonIdentitiesPractics}

In this section we briefly describe an earlier experiment to locate the $2^{27}$ periodic points of period $27$ of $p_i(z)=z^2+i$, where all but $3$ roots had been found successfully by the direct Newton Method as described above, so that the remaining $3$ roots had to be located using the Newton identities. (This is not the experiment reported earlier: subsequently we ran the experiment again with slightly improved values of $\eps_\text{stop}$ and $\eps_\text{root}$, and all roots were found; see below.) We describe how the reconstruction of the 3 missing roots among the $2^{27}$ periodic points of period $27$ of $p_i(z)=z^2+i$ works, and illustrate the results and some practical challenges. We are interested in the roots of $p(z)=p_i^{\circ n}(z)-z$ for $n=27$. The coefficients are Gaussian integers, and the first few have the following values:
\[
c_0=1, \quad c_1=0, \quad c_2= 2^{26}\,i, \quad
c_3=0, \quad c_4= -2^{51}+2^{25}+2^{25}\,i
,\quad c_5=0
\]
with $2^{51}\approx 2.25\cdot 10^{15}$; the subsequent $c_k$ with $k$ odd vanish, while those with $k$ even are huge. The respective power sums from the Newton identities have comparable sizes. Subtracting the sum of 134 million numerically computed roots to obtain power sums of the three missing roots is a certain challenge to the numerics involved.

\begin{figure}[htbp]
{
\includegraphics[width=0.95\textwidth,trim=53 585 60 72,clip]{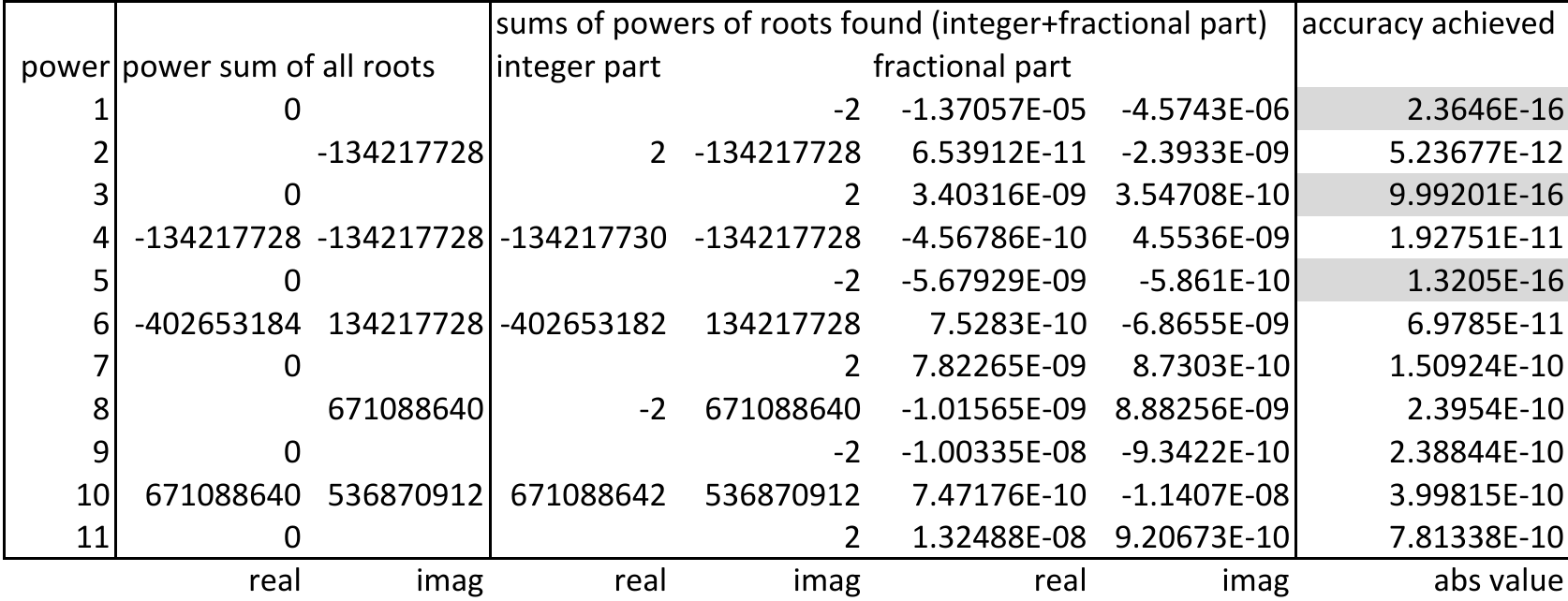}
}
\caption{Finding three missing roots using the power sum method. The first column shows the powers, the second one the power sum of all roots (using the Newton identities, separated into real and imaginary parts), the third the power sum of all $2^{27}-3$ roots found  (separated into integer and fractional parts, and displaying less accuracy than used), and the last column shows the numerical accuracy achieved for the sum of all $2^{27}$ roots, including the 3 recovered roots and the predicted value from column 2. The odd powers $1$, $3$, $5$ were used to reconstruct the three missing roots and thus show large computational accuracy. The remaining powers $7$, $9$, $11$ indicate that all roots have been found with about ten digits accuracy (probably much better for most roots). }
\label{Fig:ThreeMissingRoots}
\end{figure}

Nonetheless, in this case the three missing roots were found: two were close to $i$, and one was close to $0$. This was done initially using the power sums $1$, $2$, and $3$. However, we have $a_2=-134\,217\,728\,i$, so in the subtraction 9 valid decimal digits are lost, which is a severe limitation to accuracy. Therefore, the computations were redone using the power sums $1$, $3$, and $5$ for which no serious cancellations occur (the odd power sums of all roots should vanish). The results are shown in Figure~\ref{Fig:ThreeMissingRoots}. The third column shows that the odd power sums of the three missing roots must be equal to $\pm 2i$ with an accuracy of $10^{-5}$, and indeed the three recovered missing roots are approximately
\begin{align}
0 \quad  +& 13.7046\cdot 10^{-6} +4.57423 \cdot 10^{-6}i
\nonumber\\
i \quad -& 8.37799 \cdot 10^{-7} + 4.39432 \cdot 10^{-7}i
\label{Eq:ApproxLastThreeRoots}
\\
i \quad -& 8.36666 \cdot 10^{-7} - 4.39313\cdot 10^{-6} i \;.
\nonumber
\end{align}

The accuracy of these computations can be checked in the remaining power sums, as shown in the rightmost column in Figure~\ref{Fig:ThreeMissingRoots}.  The odd powers $1$, $3$, $5$ show only that the computations were self-consistent: the three roots were recovered so that these power sums are correct within a computational accuracy of about 16 digits.

\looseness-1
The other power sums provide independent tests. The fact that they all match to about ten valid decimal digits can be taken as experimental confirmation that all roots were found to at least this precision (but there are problems when two missing roots are very close to each other; see below). We believe that the initial $2^{27}-3$ roots were found with accuracy of at least $3\cdot 10^{-16}$ as in all other computations, so that the power sums have accuracy of around $3\cdot 10^{-16}\cdot 2^{27/2}\approx 3.5\cdot 10^{-12}$, and the accuracy of the three recovered roots should be comparable to this latter value (much lower than the remaining roots because of the accumulated error in the large sums). This is consistent with the power sum for $k=2$, and we argued above why power sums for higher powers have less accuracy even when the roots themselves have small errors. (For even powers, when the power sums were close to large Gaussian integers, we again performed the computations twice, as described earlier: once to see that the sums are close to the desired Gaussian integers, and once modulo Gaussian integers to achieve higher absolute accuracy.)

\begin{figure}[htbp]
\includegraphics[width=0.47\textwidth]{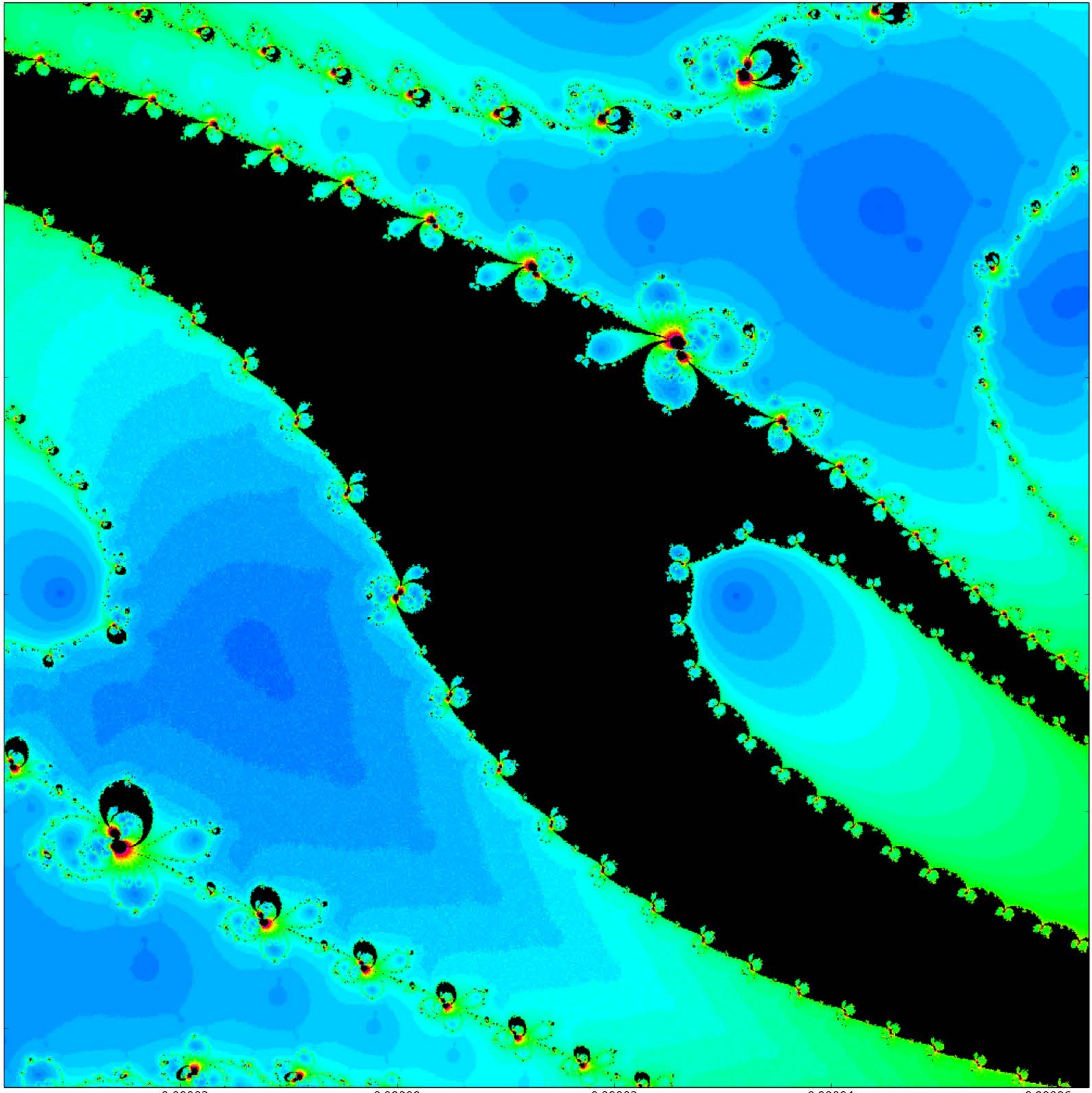} \,
\includegraphics[width=0.47\textwidth]{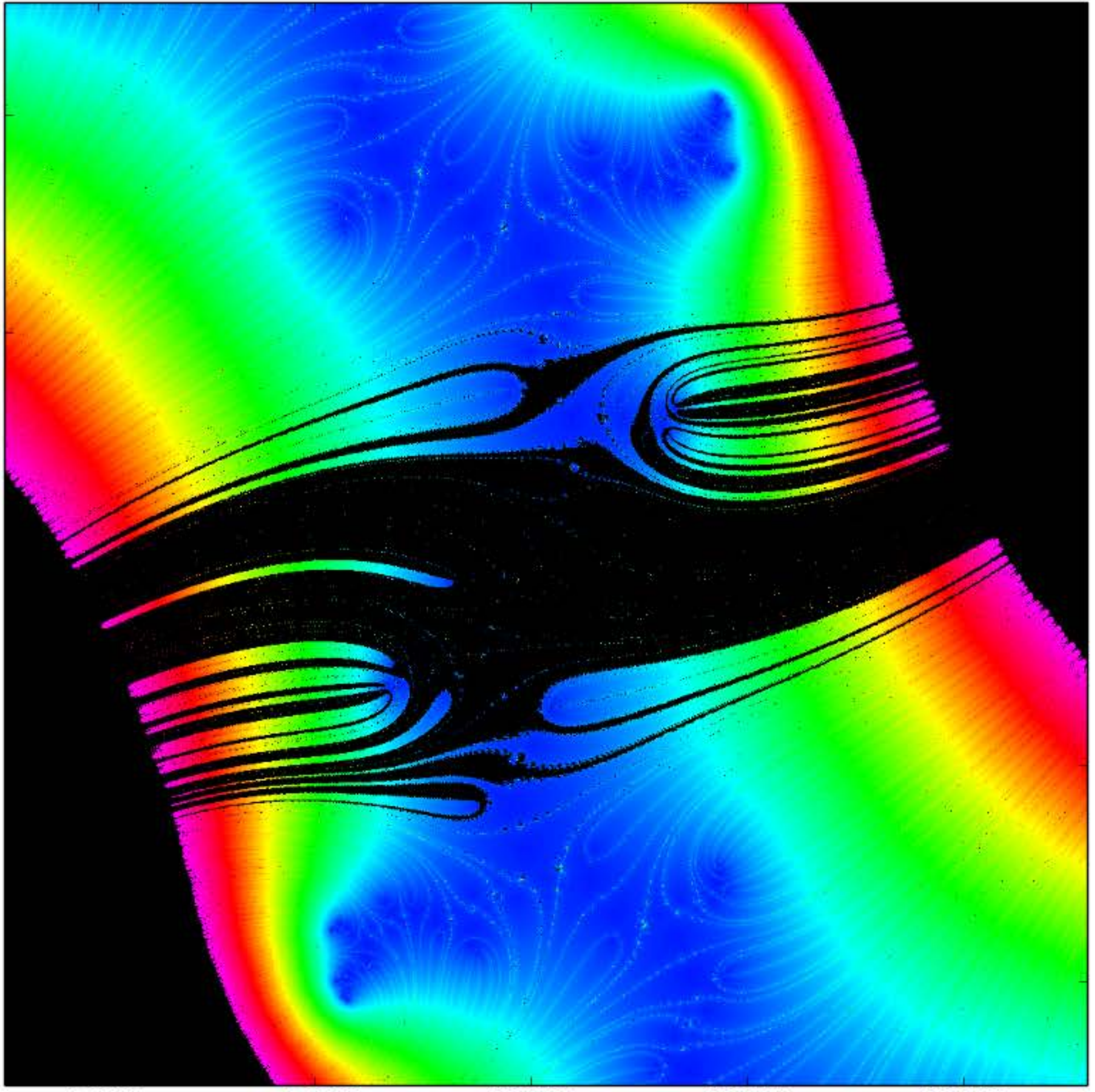}
\caption{Newton dynamics for finding periodic points of period $27$ for $z^2+i$ (degree $2^{27}$). Left: a detail near the missing root $1.37\cdot 10^{-5}+4.57\cdot 10^{-6}i$. Black shows points for which the Newton method with $\eps_\text{stop}=10^{-15}$ failed to find a root. Comparison with Figure~\ref{Fig:NewtonMarvin_1_19} (same detail for the same polynomial, with $\eps_\text{stop}=10^{-14}$) shows that the problem was a sub-optimal stopping criterion: most black points actually converge to the missing root, which was however not recognized. Right: the same polynomial with an even finer constant $\eps_\text{stop}=10^{-16}$ (zoomed out): many more roots fail to be found. } 
\label{Fig:NewtonMarvin_1_24}
\end{figure}

\begin{figure}[htbp]
\includegraphics[width=70mm,trim=10 130 20 30,clip]{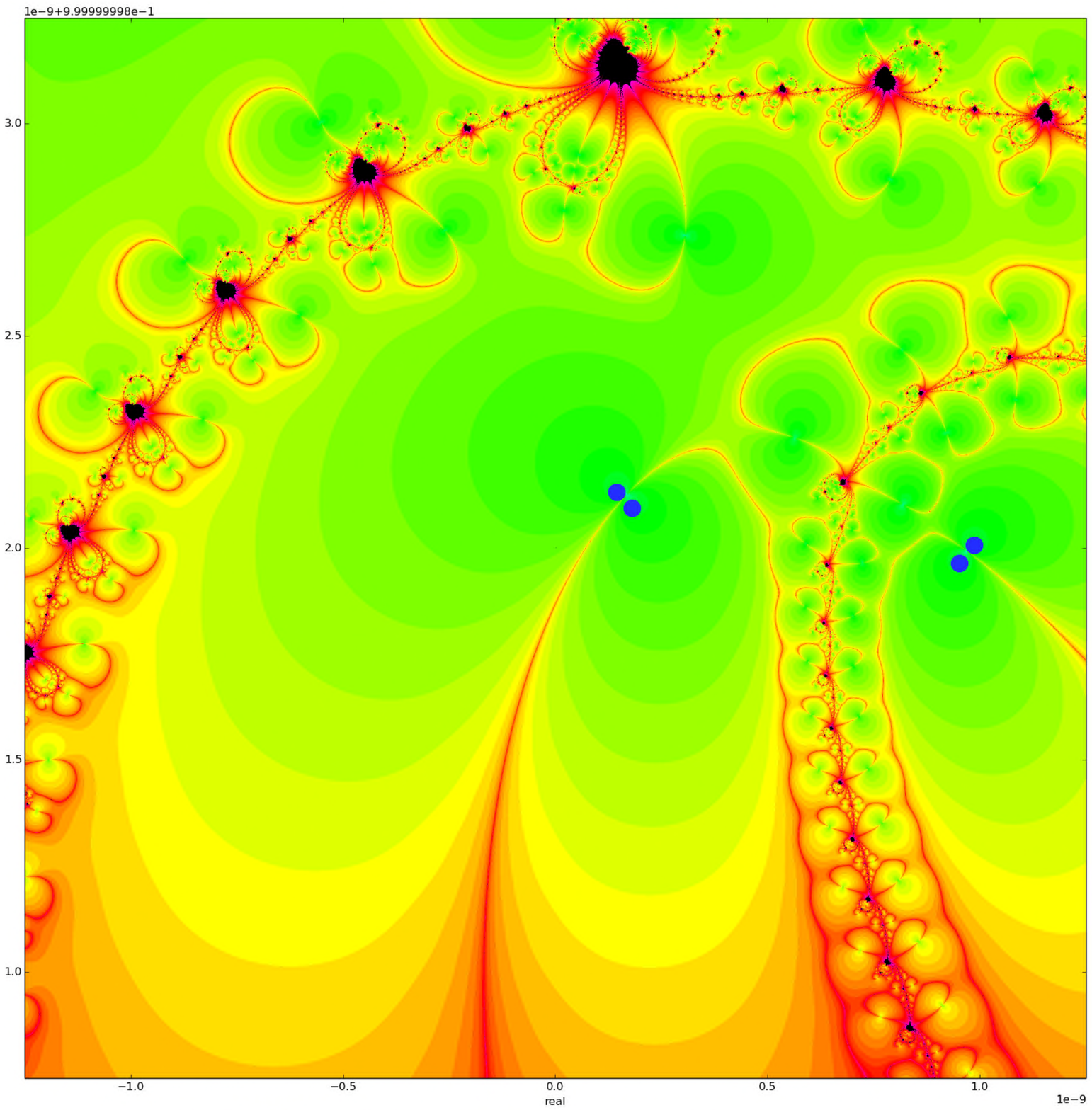}
\caption{The Newton dynamics in a neighborhood of $i$. The two missing roots near $i$ are indicated by blue dots together with the two nearby roots that were found (manually inserted, not to scale: the distance between the pairs is $8.1\cdot 10^{-10}$, while the distance within both pairs is close to $10^{-14}$). }
\label{Fig:TwoMissingRoots}
\end{figure}

With these approximations, the three missing roots were found also by a refinement of the original Newton iteration; compare Figures~\ref{Fig:NewtonMarvin_1_24} and \ref{Fig:TwoMissingRoots}. Improved approximations to these three roots are
\begin{eqnarray}
(13.704\,610\,18 +4.574\,229\,58\,i) \cdot 10^{-6} \;\;; 
\nonumber\\
i\quad +(9.676\,745\,61 - 0.072\,388\,06\,i) \cdot 10^{-10} \;;
\label{Eq:PreciseLastThreeRoots}
\\
i\quad +(1.668\,870\,81+1.253\,814\,59\,i)\cdot 10^{-10} \;\,,
\nonumber
\end{eqnarray}
now with the same experimental accuracy as all other roots.

A posteriori analysis reveals why the three roots were originally not found by our ad-hoc stopping criteria: for the first root, with our numerical accuracy we only achieve $|N(z)-z|\approx 1.33\cdot 10^{-15}>\eps_\text{stop}=10^{-15}$: we did have Newton orbits that ``wanted to'' converge to the root, but the stopping criterion was too strict. The Newton dynamics in Figures~\ref{Fig:NewtonMarvin_1_19} and \ref{Fig:NewtonMarvin_1_24} (left) illustrates this: both show the same detail of the Newton dynamics around this missing root, and both were computed with the algorithm described --- except that in Figure~\ref{Fig:NewtonMarvin_1_19} the requirement that $|N(z)-z|<\eps_\text{stop}=10^{-15}$ was relaxed to $2\cdot 10^{-15}$. The black points in Figure~\ref{Fig:NewtonMarvin_1_19} are points that do not converge to roots with the initial choice of $\eps_\text{stop}$.
Estimates using \cite[Lemma~4]{NewtonRobin1} prove that the improved approximation to the first root is indeed close to a root with error at most $7.330\cdot 10^{-15}$ (probably better).

The other two roots (both near $i$) were missed because they were closer than $\eps_\text{root}=10^{-14}$ from two other roots: in the first case, there was another root found at distance $(4.98+6.82i)\cdot 10^{-15}$, and in the second case at distance $(5.68-5.39i)\cdot 10^{-15}$, both somewhat smaller than $10^{-14}$ in absolute value.

As we pointed out, the values of the constants in our algorithm were chosen in an ad-hoc fashion and obviously have to be adjusted for polynomials of extremely large degrees --- but the underlying Newton method is stable enough that it works well even with these ad-hoc criteria. 

Once the three missing roots were found to full precision as in \eqref{Eq:PreciseLastThreeRoots}, we compared their values with the approximations from the Newton identities as in \eqref{Eq:ApproxLastThreeRoots}. For the root near $z=0$, the error was $3.50745\cdot 10^{-12}$, almost on the nose equal to the estimated value. However, for the two roots near $i$, the error turned out to be equal to 
$9.45195\cdot 10^{-7}$ for both, much larger than expected. We provide an explanation below. Of course this is below the standards of precision of the other roots found, but a good starting point for a refined search that would lead to higher precision.

Observe first that the fact that both roots have errors with almost equal absolute values is due to the fact that they must match the sums of the first powers, so the deviations must be opposite to each other within the precision of the power sums.

The fact that the last two missing roots are very close to each other has to do with a particular property of our polynomial: all the roots we are looking for are periodic points of $p_i(z)=z^2+i$, so as a set they are invariant by $p_i$. Roots near $z=0$ are contracted by the quadratic map $p_i$ that has a critical point at $z=0$ (more precisely, roots at distance $\eps$ from $0$ are mapped to roots approximately at distance $\eps^2$ from $i$). The pairs of nearest roots are thus to be expected to be close to $z=i$, and the trouble was caused by two such pairs.

After the missing root near $z=0$ is recovered numerically with expected precision, the two missing roots near $z=i$ form a quadratic equation with an (almost) double root, given with some precision $\eps$. In such a case, the roots only have precision $\sqrt\eps$, and this is what we experience here. It is an interesting fact, not without irony, that for periodic points of period $27$ in the reconstruction of the last two missing roots, an extremely ill-conditioned quadratic equation causes substantial errors, while Newton's method itself is able to find all $2^{27}$ roots (many of them much closer to each other than this pair), even the three initially missing roots and the nearby roots close to them, without any issues of ill-conditioning (except that we had to adjust our ad-hoc threshold parameters). Of course, ill-conditioning is an ubiquitous experience, but it is remarkable how well Newton can handle it.

\section{Conclusion}

The overall outcome of our three sets of experiments is that it is possible to find all roots of univariate polynomials of very large degrees in near-optimal complexity, both in terms of required Newton iterations and in terms of computing time. All roots could be located for degrees up to 1 billion ($2^{30}>10^9$), and the complexity per root is constant up to small logarithmic factors. The accuracy of all roots found is better than $3\cdot 10^{-16}$, comparable to the stopping condition $|z_n-z_{n+1}|<\eps_\text{stop}=10^{-15}$. 

In one of the intermediate experiments, for periodic points of period $n=27$ of $z^2+i$, three of the more than $134$ million roots were missed. These could be recovered with good precision by a method using the Newton identities, with remarkably little effort and good precision.

\end{document}